# Optimal sizing of renewable energy storage: A comparative study of hydrogen and battery system considering degradation and seasonal storage

Son Tay Le*, Tuan Ngoc Nguyen*,†, Dac-Khuong Bui, Tuan Duc Ngo†

*Department of Infrastructure Engineering, The University of Melbourne, Parkville, VIC 3010, Australia*

## Abstract

Renewable energy storage (RES) is essential to address the intermittence issues of renewable energy systems, thereby enhancing the system stability and reliability. This study presents an optimisation study of sizing and operational strategy parameters of a grid-connected photovoltaic (PV)-hydrogen/battery systems using a Multi-Objective Modified Firefly Algorithm (MOMFA). An operational strategy that utilises the ability of hydrogen to store energy over a long time was also investigated. The proposed method was applied to a real-world distributed energy project located in the tropical climate zone. To further demonstrate the robustness and versatility of the method, another synthetic test case was examined for a location in the subtropical weather zone, which has a high seasonal mismatch. The performance of the proposed MOMFA method is compared with the NSGA-II method, which has been widely used to design renewable energy storage systems in the literature. The result shows that MOMFA is more accurate and robust than NSGA-II owing to the complex and dynamic nature of energy storage system. Due to abundant solar resources, the system in the tropical zone always gives a superior return when compared to a similar system in the subtropical zone, highlighting the impact of weather and solar radiation. The optimisation results show that battery storage systems, as a mature technology, yield better economic performance than current hydrogen storage systems. However, it is proven that hydrogen storage systems provide better techno-economic performance and can be a viable long-term storage solution when high penetration of renewable energy is required. The study also proves that the proposed long-term operational strategy can lower component degradation, enhance efficiency, and increase the total economic performance of hydrogen storage systems. The findings of this study can support the implementation of energy storage systems for renewable energy.

---

* These authors contributed equally to this work.

† Authors to whom correspondence should be addressed. E-mail addresses: tuan.nguyen@unimelb.edu.au (Tuan Ngoc Nguyen); dtngo@unimelb.edu.au (Tuan Duc Ngo)

Keywords: Renewable energy storage; Hydrogen energy storage system; Battery energy storage system; Optimal sizing; Operation strategy; Firefly Optimisation Algorithm.

**Nomenclature**

| | |
|---|---|
| $B_y$ | electricity bill at year $y$ |
| $C_{Cap}$ | capital cost |
| $C_{O\&M,y}$ | annual O&M cost at year y |
| $C_{Rep}$ | replacement cost |
| $E_B^t$ | energy stored in the battery at time t |
| $E_{g,\text{off-peak}}^t$ | energy used from the utility grid during one-hour interval $t$ at off-peak rate |
| $E_{g,peak}^t$ | energy used from the utility grid during one-hour interval $t$ at peak rate |
| $E_{g,\text{shoulder}}^t$ | energy used from the utility grid during one-hour interval $t$ at shoulder rate |
| $I^t$ | current at time t (A) |
| $k_y$ | electricity cost factor |
| $m_{H_2}^t$ | hydrogen flow rate (kg/s) |
| $P_{charge}^t$ | charge power at time t |
| $P_{discharge}^t$ | discharge power at time t |
| $P_{EL,cell}^t$ | maximum cell power of EL at time t |
| $P_{FC,cell}^t$ | maximum cell power of FC at time t |
| $P_{G-imp}^t$ | power import from the grid at time t |
| $P_L^t$ | load at time t |
| $P_{PV}^t$ | output power from solar PV at time t |
| $r_{O\&M,i}$ | O&M cost factor for component i |
| $r_{rep,i}$ | replacement cost factor for component i |
| $r_y$ | electricity price at year y |
| $R_y$ | system revenue in year $y$ |
| $S_i$ | size of component $i$ |
| $SOH_{PEMEL}^t$ | state of health of PEMEL at time t |
| $SOH_{PEMEL}^t$ | state of health of PEMFC at time t |
| $U_i$ | unit capital cost for component i |
| $U_i$ | unit capital cost of component $i$ |
| $\eta_B^y$ | BESS efficiency at year y |
| $\eta_I$ | inverter efficiency |
| LIMIT_CLOUDY | storage limit during cloudy period |
| LIMIT_SUNNY | storage limit during sunny period |
| $n$ | number of cells |
| $t_{end}$ | end time of the sunny period |
| $t_{start}$ | start time of the sunny period |
| *MOMFA* | |
| $r_{ij}$ | distance between any two fireflies i and j |
| $x_i, x_j$ | coordinate of the $i$ and $j$ firefly |
| $X_n$ | logistic chaotic value for the $n^{th}$ firefly |
| $\alpha^t$ | trade-off coefficient at $t$ iteration |
| $\beta_0$ | firefly attractiveness at r = 0 |
| $\beta_{chaos}^t$ | $\beta_{chaos}^t$ is the chaotic number at $t^{th}$ |
| $n$ | number of the firefly |
| $t$ | number of iterations |
| $\beta$ | firefly attractiveness at each iteration |
| $\gamma$ | absorption coefficient ($0 \leq \gamma \leq 1$) |
| $\eta$ | chaotic parameter |
| *Abbreviations* | |
| BESS | Battery Energy Storage System |
| HESS | Hydrogen Energy Storage System |
| NPC | Net Present Cost |
| NPV | Net Present Value |
| PEMEL | Proton Exchange Membrane Electrolyser |
| PEMFC | Proton Exchange Membrane Fuel Cell |
| PV | Photovoltaic |
| SSR | Self Sufficiency Ratio |

# 1 Introduction

Fossil fuels provide for about 80% of the world's primary energy supply and global energy consumption is expected to increase at a rate of around 2.3% per year from 2015 to 2040 [1].

Burning fossil fuels not only threatens to increase $CO_2$ levels in the atmosphere but also emits other environmental pollutants such as $SO_x$, $NO_x$, particulate matter (PM), volatile organic compounds, and toxic heavy metals [2], escalating global warming. The increase in temperature has impacted the weather patterns and disrupted the natural balance [3]. Humans and all life forms on Earth are at risk because of this. Renewable energy sources (RESs) and renewable energy generations (REG) are believed to be alternative energy sources to fossil fuels, not only to minimise pollution and global warming but also to aid in the development of the economy. Although RESs are rapidly growing, there are challenges associated with their wider adoption of RESs. Firstly, RESs are dependent on geographic location and weather conditions. Natural resources that are used to generate renewable energy power are uncontrollable. Secondly, intermittent renewable energy can cause supply and demand mismatches, power quality issues, and network constraints [4]. Therefore, there is a need to use Energy Storage Systems (ESS) to store energy at one time and use it later. At the present, ESS is also commonly used in other applications such as demand load-shifting, PV curtailment reduction, and demand peak shaving [5].

There are several different technologies and methods to store energy. Storage capacity and discharge time are two main characteristics of energy storage technologies. Readers are encouraged to refer to previous studies [6-8] for detailed discussions on storage methods. Electro-chemical technologies allow electrical and chemical energy to be converted in a minute or shorter time frame [9]. Batteries are the most well-known electrochemical energy storage and they have been widely used in transportation, electronics, and power grid applications. As a mature technology, Battery Energy Storage Systems (BESSs) are flexible, reliable, economical, and responsive for storing energy [10-12]. Batteries, as a fast-responding ESS, are excellent at smoothing load fluctuations and absorbing short-term changes. However, as the penetration of renewable energy increases and the connectivity to the grid phases out, long-duration energy storage becomes much more important [13]. A net-zero-carbon grid needs an ESS that can store energy for long enough. Hydrogen Energy Storage System (HESS) is the system that converts energy into a form of hydrogen by either physical-based methods or material/chemical-based methods. Hydrogen has been identified as one of the key technology solutions to address climate change due to its abundant availability, high mass-energy density, and pollution-free production process. The use of hydrogen as a clean fuel as well as a long-term flexible energy storage option for backing up intermittent renewable sources has been rapidly increasing [14]. As a result, BESS and HESS as ESS options are investigated in this study.

Research on sizing EES components has been rapidly developing during the past decade. A review of BESS was accomplished by Hannan et al. [15], covering technologies, optimization objectives, constraints, approaches, and outstanding issues. Similarly, the status of research on sizing the HESS is reviewed by Xu et al. [16] and Eriksson and Gray [13]. For both ESS options, the literature review indicates that sizing optimisation is critical for the efficient and cost-effective utilisation of RES. However, most of the studies focus on optimising the size of the storage system with a predefined and fixed operational strategy. Especially, these approaches are more suitable for short term storage purposes and cannot utilise the hydrogen capability to store energy over the long term. There are only a few studies that attempt to optimize the storage capacity and the operation strategy at the same time. Sukumar et al. [17] developed a mix-mode energy management technique to lower the total cost of operation of a grid-connected BESS. A sizing method determining the optimal capacity of BESS by using the particle swarm optimisation (PSO) algorithm is also presented. In a study by Zhang et al. [18], the optimal size of grid-connected BESS was determined while optimal operations is also investigated. Pan et al. [19] proposed a bi-level mixed integer planning model that can handle the intra-day and inter-day optimization problems of HESS. They proposed a method that includes different types of decision variables, namely investment variables (binary variables that indicate the configuration of the system), operation variables (including component capacity, charge/discharge power, etc), and state variables (the ON/OFF states of system components). Their results show that the model can lower hydrogen supply prices by rationally utilising different renewable energy sources. Zhang et al. [10] use a Genetic Algorithm (GA) to simultaneously optimise the size of HESS and operational parameters. Three hydrogen storage operation techniques are examined under two scenarios based on pessimistic and optimistic costs: conventional operation strategy, peak shaving strategy, and hybrid operation strategy. The results show that the hybrid operation strategy, which combines the conventional with the peak shaving strategy, achieves the best performance.

Nevertheless, these previous studies did not consider the dynamic of other components of the system such as the variation in future electricity prices or the degradation of components throughout the analysis period. Because the performance of the battery, fuel cell, and electrolyser degrades with time, it has been shown that aging has a vital impact on system operation [15, 20, 21]. Moghaddam et al. [22] used updated forecast data to reduce the required size of an optimal BESS. In the same study, the authors also introduced a formulation to estimate the battery lifetime and used it to analyse the impact of the proposed method on the battery lifetime depreciation. Adam and Miyauchi [23] proposed different operational strategies for a hybrid battery/hydrogen energy storage system. The results show that the efficiency of peak shaving control is 52.3%, which is higher than the one obtained from the traditional power management system strategy (40.38%). The battery depreciation is also reduced from 1.2% to 1.3%. Li et al. [21] designed a combined cooling/heating/power and hydrogen microgrid system considering hydrogen degradation. Guinot et al. [24] investigate the component degradation in a hybrid battery-hydrogen system. Both size and operational strategy optimisations are carried out. From the literature, it has been shown that an optimised operation strategy can achieve superior performance when the operation strategy has different operation modes in the winter and the summer to utilise the ability of hydrogen for seasonal storage. Therefore, in this study, a model that simultaneously optimises the operational strategy and system capacity is performed considering the dynamic performance throughout the project lifetime of 25 years.

Given the complex and dynamic nature of sizing BESS/HESS components with the consideration of component degradations, price variations and operation strategies, effective optimisation algorithms are commonly required. Evolutionary algorithm-based techniques have been widely used in the design and operational control of ESS [11, 13, 25], mainly to optimise a single objective such as Net present cost [26] or Levelized cost of energy [27]. However, sizing ESS systems usually involves multi objectives (e.g., technical and economical objectives), which are mutually contradictory in most cases [14]. For instance, the increase in power reliability, which is a technical objective, will increase the system cost, which adversely impacts the economic objectives. This has led to the use of multi-objective optimisation such as Non-dominated Sorting Genetic Algorithm II (NSGA-II) [10, 16] and Multi-Objective Particle Swarm Optimization (MOPSO) [28] in sizing ESS in the literature. This study implements the Multi-Objective Modified Firefly Algorithm (MOMFA) for optimal sizing and operation of ESS. Modified Firefly Algorithm has demonstrated superior performance in optimising engineering problems in our previous work [29-31].

The aim of this study is to develop an optimisation framework for the sizing and operation of grid-connected renewable energy systems, considering the impact of component degradation, price variations and operation strategies using MOMFA. Two storage technologies, namely BESS and HESS are considered for a case study of a warehouse using real monitored electricity usage data and measured solar irradiation. Major contributions of this paper are:

- The proposed framework optimises the size and energy management of a renewable energy system with a consideration of electricity price forecasting, solar output depreciation, and degradation of storage components in the system.
- A new operational strategy with two distinct operation modes in the summer and winter is introduced to employ the potential of hydrogen to store energy over a long period.
- A new multi-objective Modified Firefly Algorithm (MOMFA) is applied to obtained Pareto set.
- The proposed approach is applied in a real-world case study of a warehouse in the tropical zone, which has year-round sunshine. Another test case is examined for the subtropical weather zone with a high seasonal mismatch to further verify the robustness and versatility of the method.
- The economic performance of HESS is compared with BESS, which is a mature and widespread used energy storage technology, to evaluate its market readiness levels and future potential.

The paper is structured as follows. The methodology is described in Section 2. The proposed framework is applied to a real-world case study in Section 3. The result and discussion are reported in Section 4. Section 5 concludes the findings of the study.

# 2 Methodology

This section describes the system configurations as well as the proposed optimal sizing method.

## 2.1 System configurations

Three system configurations are investigated based on different storage technologies and operational strategies. Two proposed grid-connected renewable energy systems are shown in Figure 1. The main components are a rooftop solar PV system, electricity load, an external grid and an energy storage system (BESS in Figure 1a or HESS in Figure 1b). In Case 1 and 2, the electricity management of both ESS follows a conventional strategy (CS) which is widely adopted in the literature [10, 16, 26, 32]. In addition, to assess the potential of HESS for long-term storage, an optimised long duration strategy (OLDS) is proposed in Option 3. The details and operational strategies of BESS and HESS will be discussed in the next section.

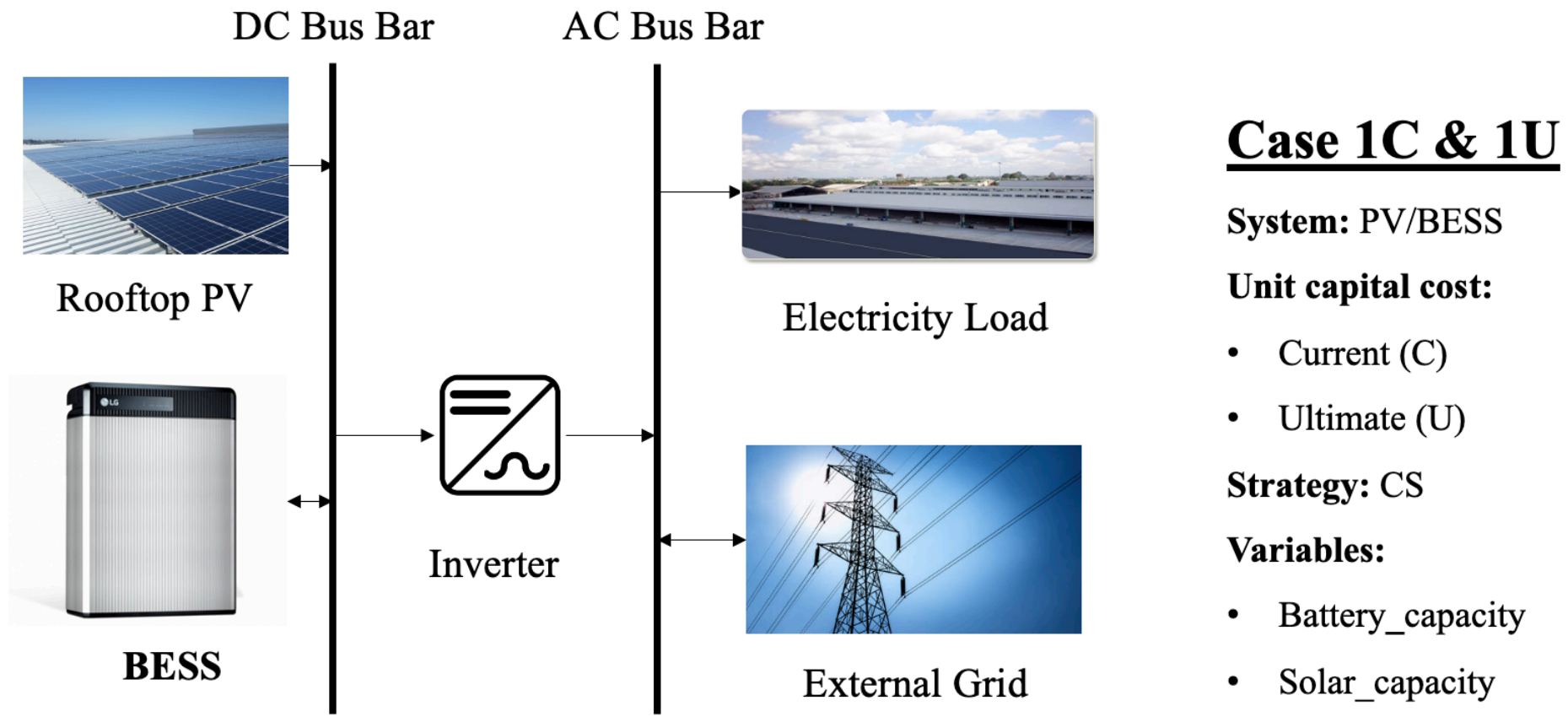

(a) PV/BESS

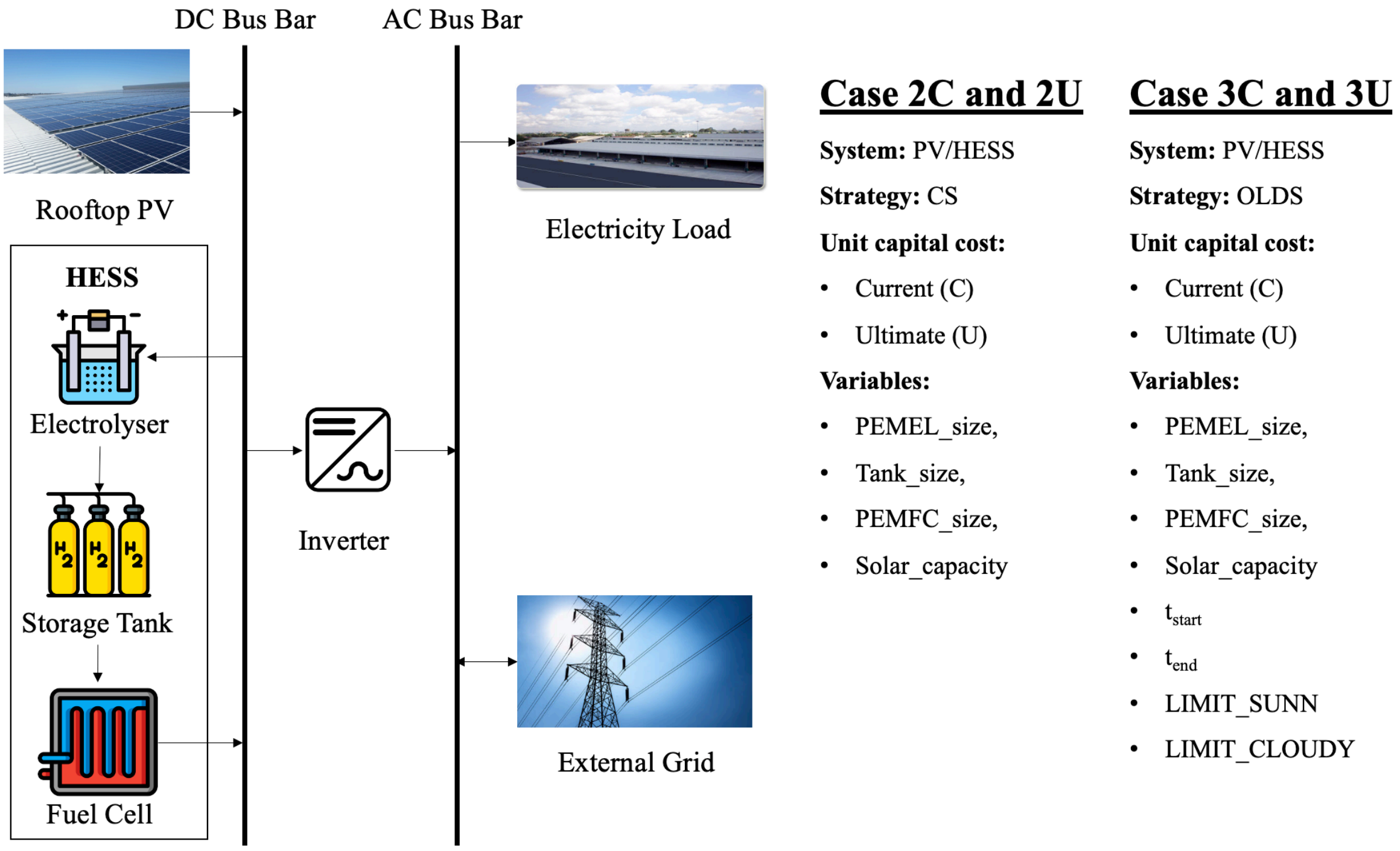

(b) PV/HESS

Figure 1 System schematic layout of the two proposed grid-connected systems

Table 1 summarises the cost of energy storage components from the literature review. As observed in Table 1, the cost of BESS and HESS varies considerably. The cost of each component is greatly influenced not only by technology and overall capacity but also by the country of the manufacturers. As the result, two cost scenarios are used in this study, namely current cost (C) and ultimate cost (U). The current cost scenario is based on the previous review of the component cost. On the other hand, the ultimate cost scenario is based on the cost target set by the International Energy Agency (IEA) [33] and the United States Department of Energy (DOE) [34, 35].

Table 1 Review Unit Cost

| Component | CAPEX Price Range | O&M Price Range | Ref: Location/Lab – Year – Cost |
|---|---|---|---|
| Lithium-ion battery | 200 USD/kWh – 800 USD/kWh: depends on technology and total capacity. | 0.5~3% CAPEX | • United States – 2017 - 200~500 USD/kWh [36]<br>• Spain – 2018 – 800 USD/kWh [37]<br>• Italy – 2021 – 640 USD/kWh [38]<br>• UK – 2021 - 300~380 USD/kWh [39]<br>• Sweden – 2017 - 400 USD/kWh [10] |
| Electrolyser | 300 USD/kW – 5,000 USD/kW: depends on technology and total capacity. | 1~5% CAPEX | • Qatar - 2020 - 300 USD/kW [40]<br>• DoE – 2021 - 460 USD /kW [41]<br>• IRENA – 2020 - 650~1,000 USD/kW [42]<br>• Brazil - 2018 – 1,500 USD/kW [43]<br>• Australia – 2018 - 2,000 USD/ kW [44]<br>• Sweden – 2017 - 5,000 USD/ kW [10]<br>• EU – 2018 - 985~2,250 USD/kW [45] |
| Storage tank | 300 USD/kg - 600 USD/kg | 0.5 – 2.5% CAPEX | • Qatar - 2020 - 300 USD/kg [40]<br>• DoE – 2021 - 378 USD /kg [41]<br>• Brazil - 2018 - 500 USD/kg [43]<br>• Australia – 2018 - 438 USD/ kW [44]<br>• Sweden – 2017 - 576 USD/ kW [10] |
| Fuel cell | 100 USD/kW - 5,600 USD/kW: also depends on the capacity. The cost of fuel cell is one of the | 1~5% CAPEX | • Qatar - 2020 - 100 USD/kW [40]<br>• DoE – 2021 - 193 USD /kW [41]<br>• Brazil - 2018 - 2,000 USD/kW [43]<br>• Australia – 2018 - 600 USD/ kW [44]<br>• Sweden – 2017 - 4,000 USD/ kW [10] |

| | | |
|---|---|---|
| | most significant challenging in commercializing HESSs. | • India - 2016 - 1,600 USD /kW [46]<br>• Australia – 2021 - 5,600 USD/kW [47] |

### 2.1.1 Battery Energy Storage System

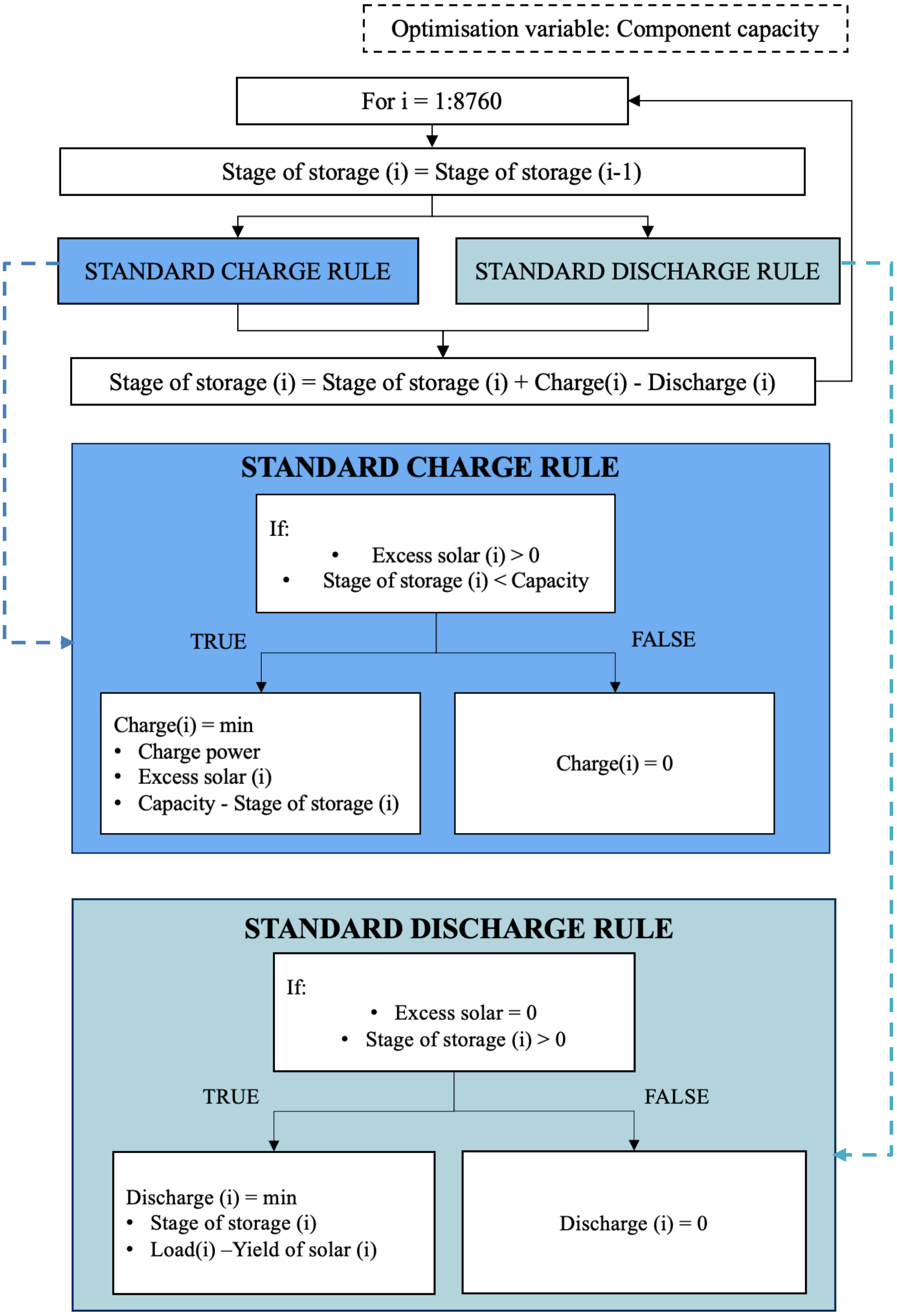

Figure 2 Flowchart of conventional strategy (CS)

BESS has been widely used to store excess energy from renewable energy systems. As a result, the efficiency of BESSs has been tested extensively [48]. In this paper, a BESS with constant annual efficiencies is modelled based on the conventional strategy (CS), which has been implemented in many studies [27, 49, 50]. The flowchart of the CS is presented in Figure 2, which shows that as long as the stage of the storage is less than the total capacity, the storage system will be charged if there is excess solar. The energy stored in the battery ($E_B$) may be calculated using the following equation:

$$E_B^t = E_B^{t-1} + P_{charge}^t \times \Delta t \tag{1}$$

$$P_{charge}^t = \max\left(P_{PV}^t - \frac{P_L^t}{\eta_I}, 0\right) \tag{2}$$

Under the discharging stage, the storage discharges when the power from the solar PV is not enough to supply the electricity demand, given that the stage of storage is greater than zero. The energy available in the battery is updated using the following equations:

$$E_B^t = E_B^{t-1} - P_{discharge}^t \times \Delta t \tag{3}$$

$$P_{discharge}^t = \max\left(\frac{1}{\eta_B^y}\left(\frac{P_L^t}{\eta_I} - P_{PV}^t\right), 0\right) \tag{4}$$

In this study, incremental time interval $\Delta t$ is 1 hour, inverter efficiency $\eta_I = 97\%$. Based on the experimental result from [48], the initial efficiency of the battery ($\eta_B^0$) is 95% and then it depreciates at a constant rate of 2.9% per year, assuming 300 cycles/year.

Table 2 summarises the cost and specification of BESS. It should be noted that in the ultimate cost scenarios, the replacement cost factor $r_{rep}$ is not used due to the limits of raw material costs.

Table 2 Battery Cost and Specification

| | **Battery** |
|---|---|
| Energy/Power (E/P) ratio | 2.5 |
| Lifetime (years) | 12 |
| Unit capital cost, $\boldsymbol{U_i}$ (USD/kW) | |
| • Current | • $490 |
| • Ultimate | • $270 |
| O&M factor, $\boldsymbol{r_{O\&M,i}}$ | 0.5% |
| Replacement cost factor, $\boldsymbol{r_{rep,i}}$ (for current scenario only) | 50% |

### 2.1.2 Hydrogen Energy Storage System

As batterie storage systems are not suitable for long-duration storage, especially for some regions such as central Europe where the storage capacity must be raised to several weeks [51]. Compared to batteries, hydrogen outperforms as a long-term storage solution as it has a higher energy density and has a low self-discharge (leakage) rate [10, 52]. Self-discharge in batteries is caused by internal chemical reactions, which reduces the energy stored without any connection between the electrodes or any external circuit. While lithium-ion batteries are found to have a self-discharge rate of about 2% to 3% per month [53], the leakage rate of hydrogen storage vessels is required to be less than 0.12% per month [54]. Therefore, the HESS option is also investigated in this study.

The HESS consists of Proton Exchange Membrane Electrolyser (PEMEL), hydrogen tank and Proton Exchange Membrane Fuel Cell (PEMFC), as shown in Figure 1b. The HESS is flexible to combine different charge power, discharge power and storage capacity because of the modularity and independence of each component. The electrolyser sets the charge power, the storage tank sets the storage capacity, and the fuel cell independently sets the discharge power. While lithium batteries have the Energy/Power (E/P) ratio is usually under eight hours [55], the E/P ratio of HESS is not constrained and may be flexibly changed depending on applications in a variety of settings. Due to the complexity of the decoupling ability of the system, each component of a HESS must be modelled in detail.

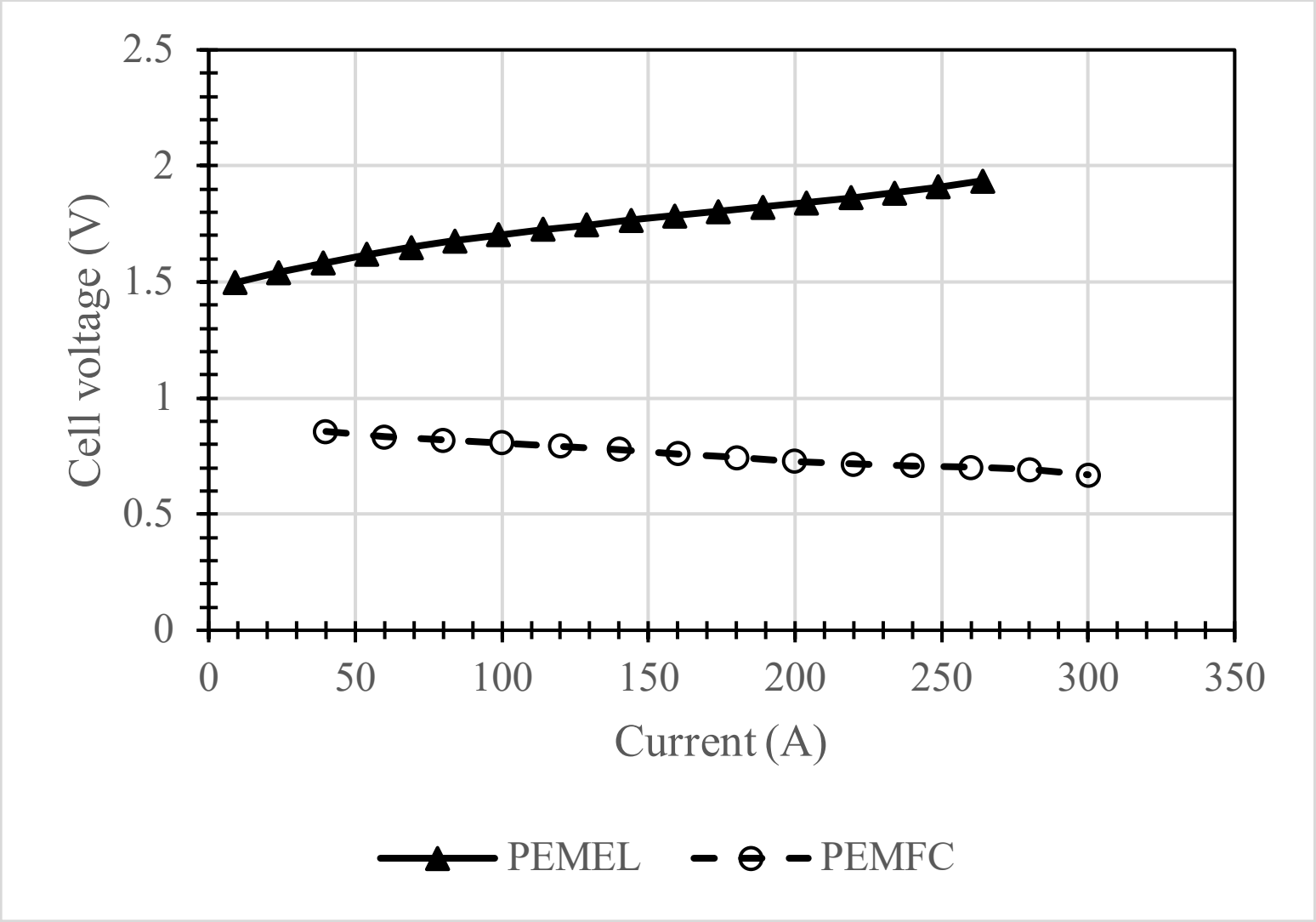

Figure 3 PEMEL and PEMFC polarization curve, derived from Guinot et al. [24]

Commercial PEMEL and PEMFC, unlike BESS, have no published testing results on their performance over time. Thus, the modelling of PEMEL and PEMFC in this paper, including the degradation model, is based on a theoretical model by Guinot et al. [24], which is based on polarization curves in Figure 3. Readers are encouraged to refer to Guinot et al. [24] for further details, including the PEMEL and PEMFC technical parameters, and comparisons between models with and without aging effects. PEMEL and PEMFC degradation model assumes a proportional increase in cell voltage with operating time. The change in cell voltage over time for PEMEL and PEMFC are $+10\mu V/h$ and $-5\mu V/h$, respectively. The state of health (SOH) of PEMEL and PEMFC is determined as:

$$SOH_{PEMEL}^{t} = \frac{\max(P_{EL,cell}^{0})}{\max(P_{EL,cell}^{t})} \tag{5}$$

$$SOH_{PEMFC}^{t} = \frac{\max(P_{FC,cell}^{t})}{\max(P_{FC,cell}^{0})} \tag{6}$$

Where $P_{EL,cell}^{0}, P_{FC,cell}^{0}$ is the initial maximum cell power of EL and FC, respectively.

As for the hydrogen tank, Eq. (7) describes the relationship between current and hydrogen flow rate [56]:

$$m_{H_2}^{t} = 1.05 \times 10^{-8}(I^{t} \times n) \tag{7}$$

The cost and specification of HESS are summarised in Table 3. The lifetimes of the electrolyser, hydrogen tank, and fuel cell are 15, 25, and 5, respectively. Therefore, during the 25-year analysis period, the electrolyser needs to be replaced once in year 15 with $r_{rep,EL} = 60\%$. Hydrogen tanks do not need to be replaced entire the simulation period and fuel cells must be replaced 4 times at a rate $r_{rep,FC} = 77.5\%,\ 55\%,\ 32.5\%,\ 10\%$ at year 5, 10, 15, 20, respectively. In ultimate cost scenarios, similar with BESS, $r_{rep} = 1$ for all components.

Table 3 HESS Component Cost and Specification

| **Component** | **PEMEL** | **Hydrogen Tank** | **PEMFC** |
|---|---|---|---|
| Lifetime (years) | 15 | 25 | 5 |
| Unit capital cost, $\boldsymbol{U}$ | USD/kW | USD/kg | USD/kW |
| • Current | • 1500 | • 600 | • 4000 |
| • Ultimate | • 200 | • 266 | • 400 |
| O&M factor, $\boldsymbol{r_{O\&M}}$ | 1% | 1% | 1% |
| Replacement cost factor, $\boldsymbol{r_{rep}}$ (for current scenario only) | 60% | - | [77.5%, 55%, 32.5%, 10%] |

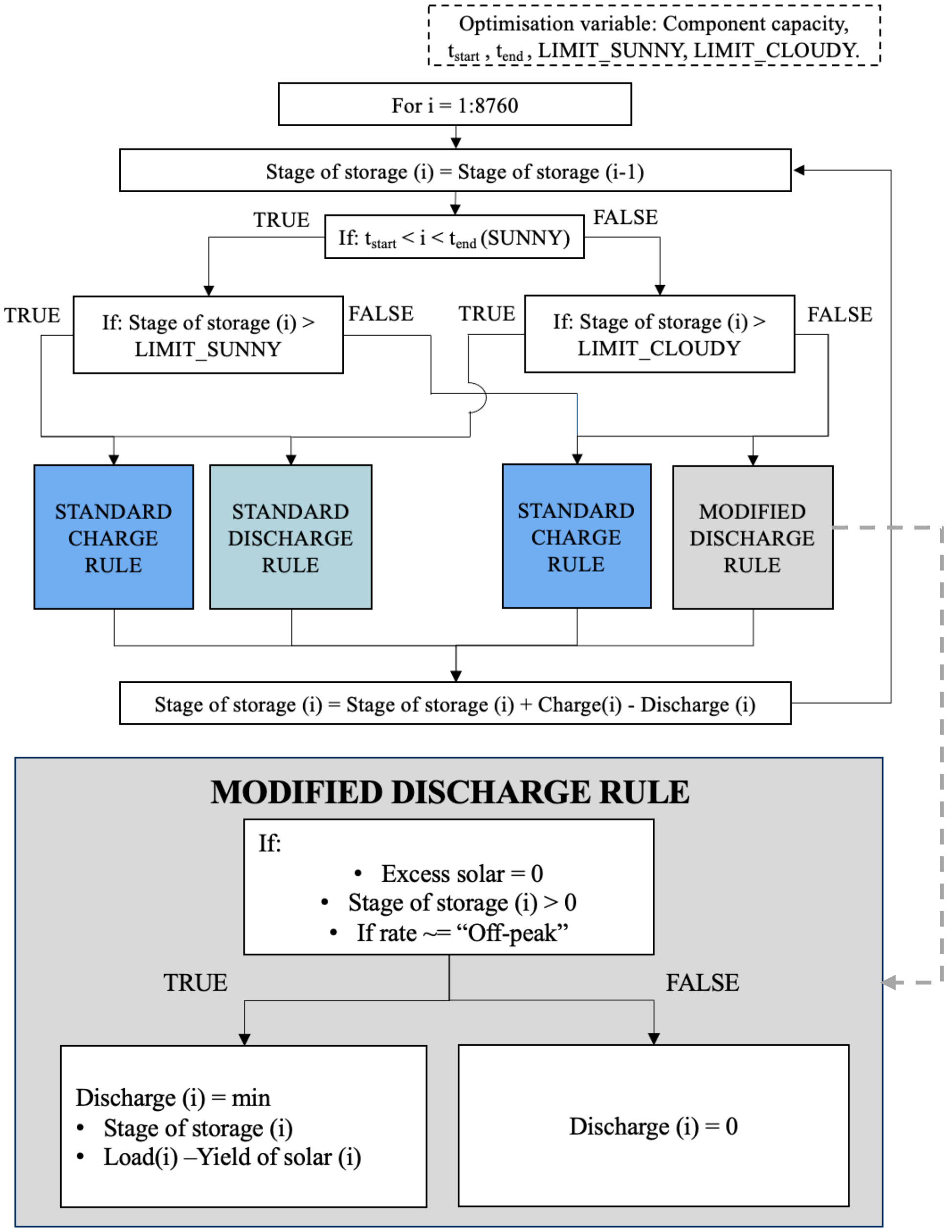

Figure 4 Flowchart of optimised long duration strategy (OLDS), refer to Figure 2 for the standard charge and discharge rule

The conventional strategy presented in the previous section for the BESS is also used to examine the performance of the HESS. In addition, the optimised long duration strategy (OLDS) is also proposed to investigate the potential of HESS for long-term storage. For the OLDS, it is assumed that hydrogen is stored during sunny months and consumed during cloudy months, utilising the long-term storage ability of HESS for seasonal storage. The flowchart of the OLDS is presented in Figure 4. As can be seen in Figure 4, two arbitrary time indicators, $t_{start}$ and $t_{end}$, are added to define the start and end of the sunny period. The idea of specifying distinct operation conditions based on $t_{start}$ and $t_{end}$ has proven to be beneficial for seasonal storage [10]. During the sunny period ($t_{start} < i < t_{end}$), if the stage of storage is greater than a defined limit (LIMIT_SUNNY), the system follows the standard charge and discharge rule similar to CS. On the other hand, if the stage of storage is lower than that limit, the system still follows a standard charge rule but uses a modified discharge rule. In this paper, the electricity supplier sells electricity at different rates, namely peak, shoulder and off-peak rates. As shown in Figure 4, in a Modified Discharge Rule, the system is similar to the standard one, except it does not discharge during off-peak hours. As the electricity price is usually low during off-peak hours, this modified rule aims to maximize the user's economic benefits while minimizing drawing power from the storage system. Similarly, another decision variable LIMIT_CLOUDY is introduced with the same principle as LIMIT_SUNNY to control the usage during cloudy months. Therefore, in the OLDS, the operational variables ($t_{start}$, $t_{end}$, LIMIT_SUNNY, LIMIT_CLOUDY), as well as the component sizing variables (e.g., PEMEL_size, Tank_size, PEMFC_size, Solar_capacity), will be optimised using the proposed optimisation algorithm presented in the following section. It should be noticed that for the conventional strategy, only the component sizing variables will be optimised.

## 2.2 Optimization method and implementation

### 2.2.1 Objective function

This study aims to optimise the component sizing variables and operation variables of the BESS/HESS as presented in the previous section to maximize the Net Present Value (NPV) and Self Sufficiency Ratio (SSR). While the NPV presents the economic incentives of the renewable energy system, the SSR can be considered as a technical-environmental goal as it indicates the percentage of electricity demands that are met by the renewable system. The computation of NPV and SSR is presented in the following section.

### 2.2.2 Net Present Value

The electricity bill for any year $y$ is calculated as:

$$B_y = \sum_{t=1}^{8760} \left(E_{g,peak}^{t} \times r_y^{peak} + E_{g,\text{shoulder}}^{t} \times r_y^{shoulder} + E_{g,\text{off-peak}}^{t} \times r_y^{off-peak}\right) \tag{8}$$

$$r_y = k_y \times r_1 \tag{9}$$

The system revenue in year $y$ is the total saving in the electricity bill due to the use of the RES:

$$R_y = B_{y,0} - B_{y,ESS} \tag{10}$$

Where $B_{y,ESS}$ is the electricity bill of the system with the expanded solar PV and ESS, and $B_{y,0}$ is the electricity bill of the system without any solar PV and ESS. In other words, $B_{y,0}$ is the electricity bill if the grid supplies 100% load. Therefore, the difference between $B_{y,0}$ and $B_{y,ESS}$ is the reduction of electricity bills when using a renewable energy system.

For both ESS and PV system, the cost consists of three costs: capital cost, O&M cost, and replacement cost. The capital cost ($C_{Cap}$) of the system, annual O&M cost ($C_{O\&M,y}$) at year y, the replacement cost $C_{Rep}$ is obtained with E.q. (11) to (13):

$$C_{Cap} = \sum_{i=1}^{n} U_i \times S_i \tag{11}$$

$$C_{O\&M,y} = \sum_{i=1}^{n} U_i \times S_i \times r_{O\&M,i} \tag{12}$$

$$C_{Rep} = \sum_{i=1}^{n} U_i \times S_i \times r_{rep,i} \tag{13}$$

The NPC and NPV over the project period of 25 years:

$$NPC = \sum_{y=1}^{25} \frac{C_{O\&M,y} + C_{Rep,i}}{(1+i)^y} + C_{Cap} \tag{14}$$

$$NPV = \sum_{y=1}^{25} \frac{R_y - C_{O\&M,y} - C_{Rep,i}}{(1+i)^y} - C_{Cap} \tag{15}$$

Where $i = 5\%$ is the discount rate.

### 2.2.3 Self Sufficiency Ratio

The SSR is the percentage of the electricity demand that is supplied by renewable energy, demonstrating the amount of renewable energy penetration:

$$SSR = \left(1 - \frac{\sum_{t=0}^{t=N} P_{G-imp}^{t}}{\sum_{t=0}^{t=N} P_{L}^{t}}\right) \times 100\% \tag{16}$$

## 2.3 Optimisation Algorithm

The operational and component sizing variables presented in the previous section will be determined using the proposed Multi-Objective Modified Firefly Algorithm (MOMFA) based on our previous work [29-31].

The conventional FA, which is a stochastic, nature-inspired metaheuristic algorithm based on the flashing characteristics and behaviors of fireflies, was firstly introduced by Yang [57]. To improve the performance of the conventional FA, Chou and Ngo [58] and Bui et al. [29] developed the modified firefly algorithm (MFA) with auxiliary components such as chaotic maps, AIW, and Lévy flight. In the MFA, the objective function of the optimisation problem is determined by the brightness of a firefly. A brighter firefly will attract other fireflies in the surrounding domain. To find the optimal solution in the MFA, the initial population is firstly generated using a Logistic map [58, 59]:

$$X_{n+1} = \eta X_n(1 - X_n), 0 \leq X_0 \leq 1 \tag{17}$$

The attractiveness of each firefly can be defined through E.q (18) to (20):

$$\beta = (\beta_{chaos}^t - \beta_0)e^{-\gamma r_{ij}^2} + \beta_0 \tag{18}$$

$$\beta_{chaos}^t = \begin{cases} 0 & \beta_{chaos}^{t-1} = 0 \\ 1 / \beta_{chaos}^{t-1} mod(1) & \text{otherwise} \end{cases} \tag{19}$$

$$r_{ij} = \left\| x_i - x_j \right\| \tag{20}$$

In this study, the chaotic parameter $\eta$ is set to 4 and absorption coefficient $\gamma$ is set to 1 based on sensitivity analysis. E.q (22) depicts the movement of firefly $i$ when it is attracted to another more attractive firefly $j$:

$$x_i^{t+1} = x_i^t + \beta\left(x_i^t - x_j^t\right) + \alpha^t sign[rand - 0.5] \otimes L\acute{e}vy \tag{22}$$

Where the second term is the attraction term which has been discussed previously; the third term is randomisation term, rand $\in [0,1]$ is a random number generated by a uniform distribution in [0, 1] and $\otimes$ is entry-wise multiplication. $\alpha^t$ can be calculated with E.q. (23).

$$\alpha^t = \alpha_0 \theta^t \tag{23}$$

where $\alpha_0 = 1$ is the initial randomisation parameter; $\theta$=0.9 is the randomness reduction constant based on sensitivity analysis and the literature.

Lévy flights are determined as follows:

$$L\acute{e}vy \sim s = \frac{u}{[v]^{1/\tau}} \tag{24}$$

where $u$ and $v$ are determined by a normal distribution:

$$v \sim N(0,1) \tag{25}$$

$$u \sim N(0, \left\{ \frac{\Gamma(1+\tau)\sin\left(\frac{\pi\tau}{2}\right)}{\Gamma\left[\frac{(1+\tau)}{2}\right]\tau 2^{\frac{\tau-1}{2}}} \right\}^{2/\tau}) \tag{26}$$

where is $\tau = 3/2$; $\Gamma$ is the Gamma function as determined by E.q. (27):

$$\Gamma(z) = \int_0^\infty t^{z-1} e^{-t} dt \tag{27}$$

For the multi-objective optimisation problem in this study (i.e. maximising NPV and SSR), a Pareto optimal solution is implemented. The MOMFA will perform the optimal search to find solutions in the Pareto optimal set. A multi objective's Pareto Front (PF) can be defined as the set of non-dominated solutions:

$$PF = \{s \in S \mid \exists s' \in S : s' \prec s\} \tag{28}$$

Where $S$ is the solution set, $\prec$ is the non-dominance relationship. Based on Pareto dominance, solutions in a population P are separated into disjoint subsets and ranked based on a Non-Dominated Sorting Algorithm (NDSA). Readers are encouraged to refer to [60] for a detailed description of NDSA.

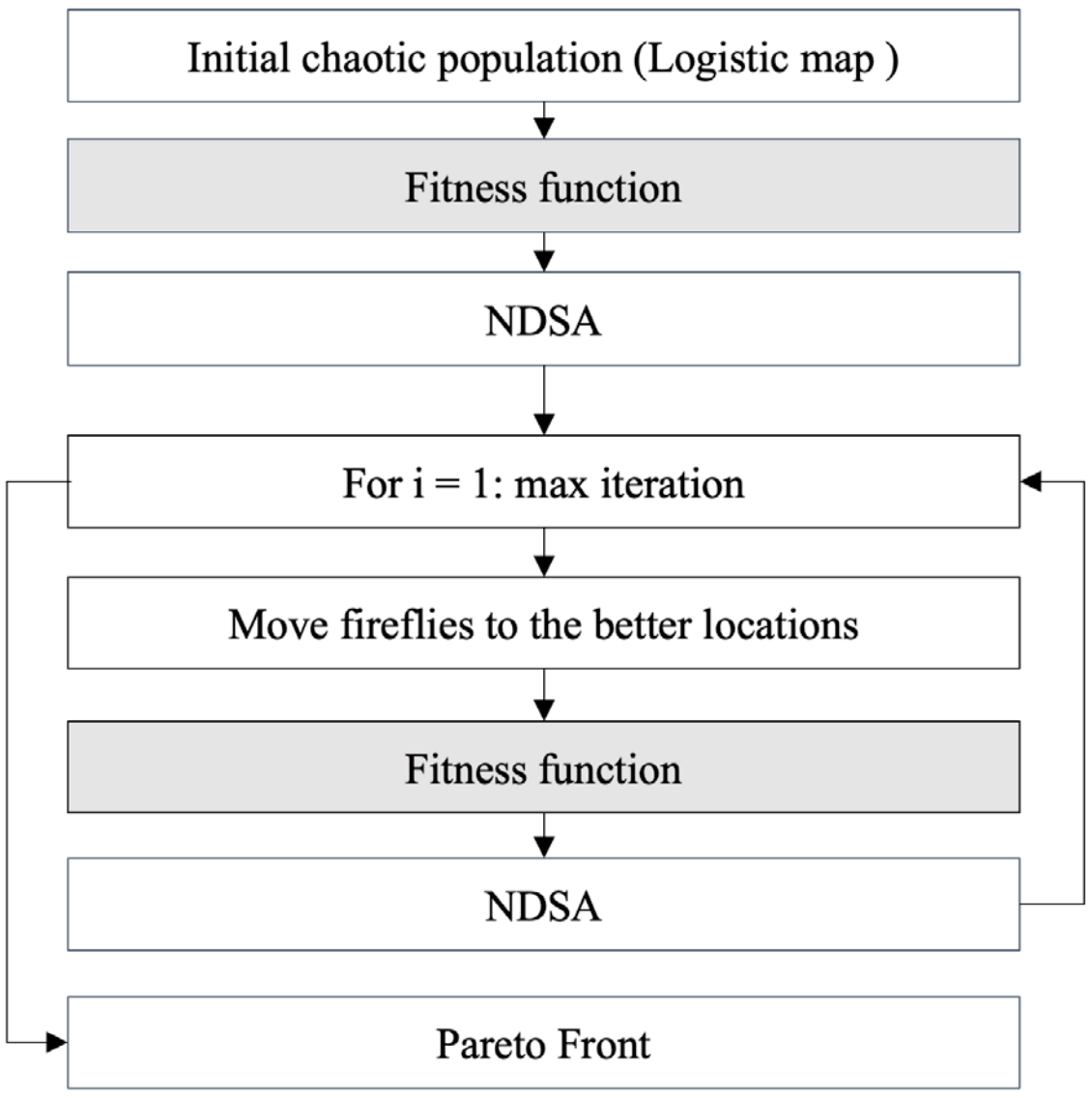

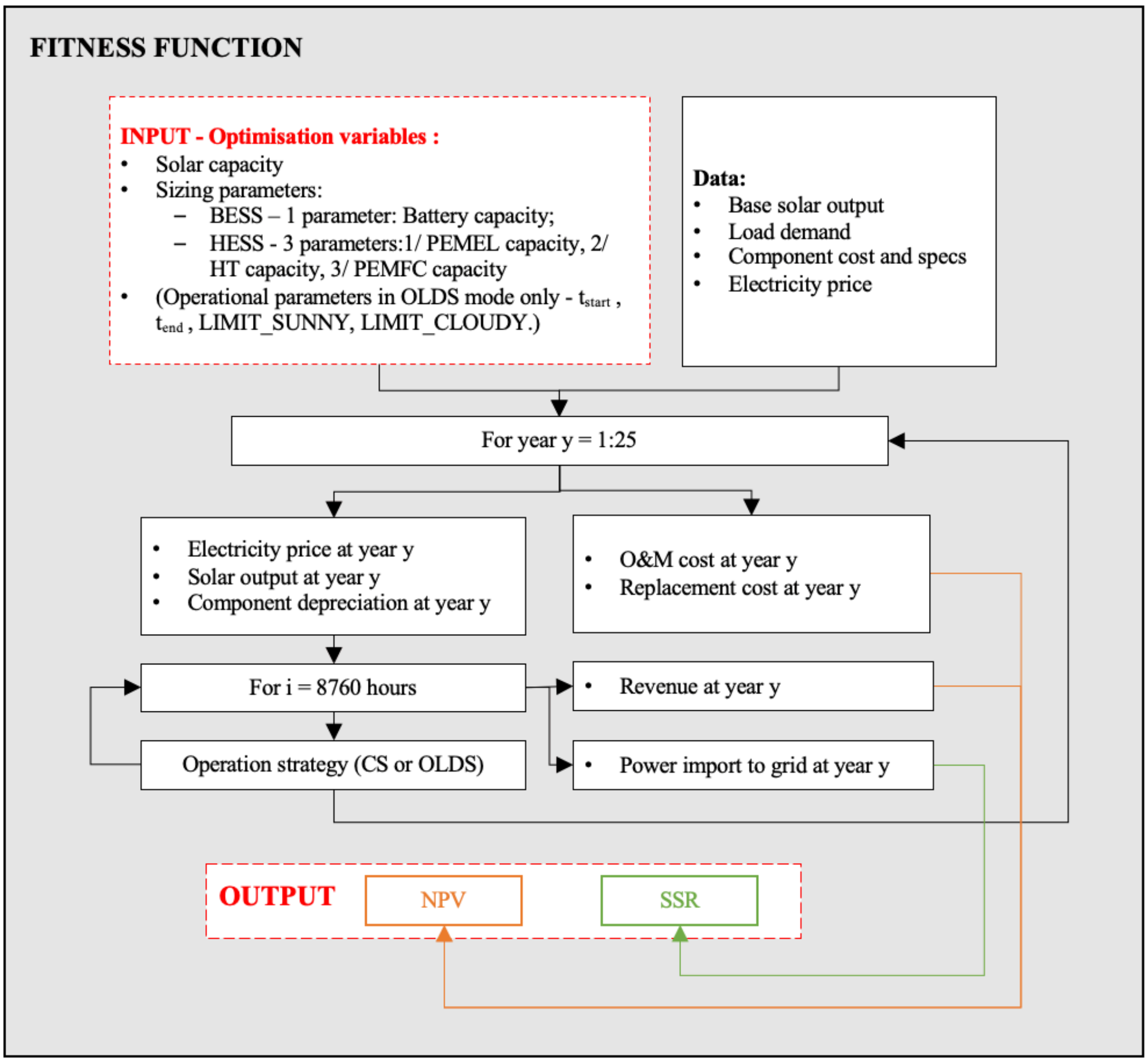

Figure 5 MOMFA Flow Chart

Figure 5 depicts the entire flowchart of the MOMFA, which is implemented in MATLAB 2021b for optimising the sizing and operational variables of BESS/HESS in this study. The optimal results obtained from the MOMFA will be compared with those obtained from the NSGA-II code [61].

# 3 Case study

The case study is a distributed energy project in Ho Chi Minh (HCM) city, Vietnam. The warehouse has a rooftop PV system, which provides 1,5 million kWh/year, as a way of reducing energy costs while also fostering environmental stewardship. In this study, the PV system will be expanded, and excess energy will be stored in an energy storage system. The study aims to propose the optimal sizing of the PV-energy storage system through the use of an optimisation design to increase the penetration and efficiency of renewable energy systems. It is noted that solar rebates and feed-in tariffs are not considered in this study.

The PV and electricity usage in this study is measured with real monitored data. Figure 6 presents the actual load consumption of the warehouse for a typical week. The daily yield of the current PV system is also presented in Figure 6. As aforementioned, this study investigates the potential of expanding PV and installing storage systems for the higher penetration and efficiency of renewable energy. The yield of expanded PV systems is assumed to be linearly proportional to the base solar. In addition, a linear depreciation rate of 0.55% per year for the output of the PV system is also applied in this study. The unit cost of the expanded solar system is \$881/kW, which is based on the total installed cost of commercial sector solar PV in South East Asia [62]. The lifetime of the solar PV is 25 years, and the O&M cost of the solar system is assumed to be 1% of CAPEX per year.

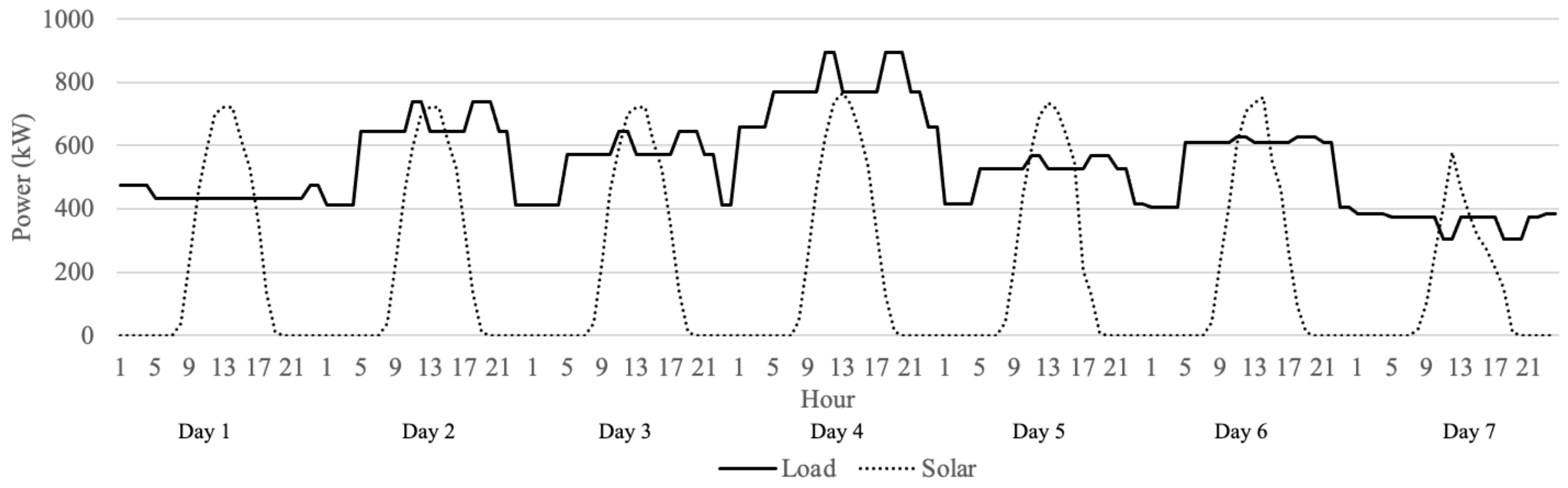

Figure 6 Sample hourly load profile and current solar PV yield power

Moreover, the warehouse is located in HCM which has a tropical climate, and the solar radiation varies little throughout the year. This study also investigates a hypothetical situation, where the warehouse is in Melbourne, Australia (MEL) - a location with a high seasonal mismatch problem. All input data for the optimal design of BESS and HESS remain the same except for the change in solar power output. HOMER PRO Software [63] is used to estimate the solar power output in MEL. Readers are encouraged to refer to [64] for more details on HOMER software functionalities and implementations. Figure 7 presents the comparison of real solar yield in HCM and HOMER-generated solar yield in MEL.

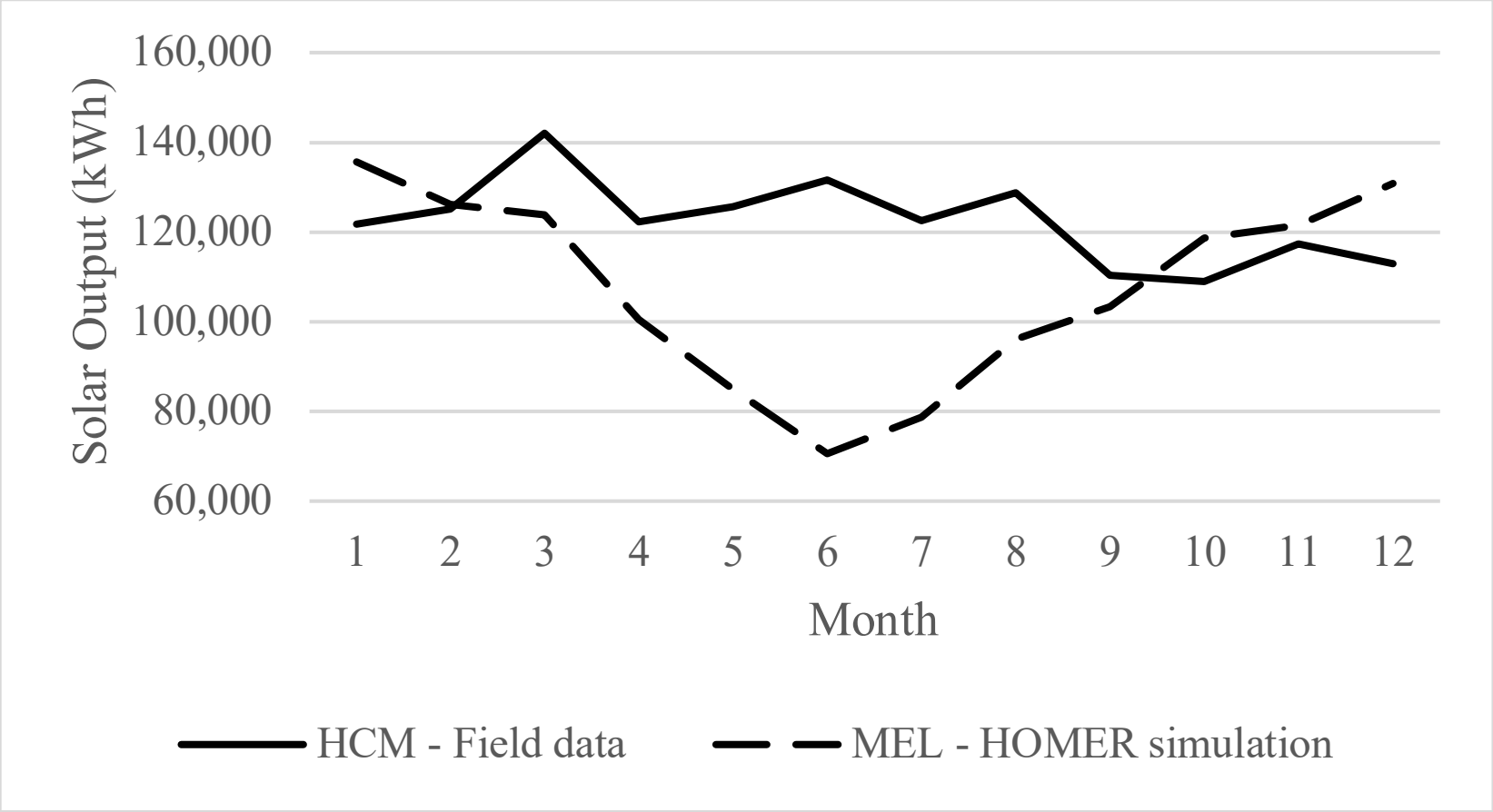

Figure 7 Monthly solar output in HCM and MEL

Table 4 shows the definition and the corresponding electricity price of each rate in the first year and the electricity cost factor $k_y$ is shown in Figure 8.

Table 4 Electricity Rate

| **Rate** | **Period** | **Electricity Price (USD/kWh)** |
|---|---|---|
| Peak | **Mon-Sat:** 10:00 to 12:00 (2 hours); 17:00 to 20:00 (3 hours) | $r_1^{peak} = 0.187$ |
| Shoulder | **Mon-Sat:** 4:00 to 10:00 (6 hours); 12:00 to 17:00 (5 hours); 20:00 to 22:00 (2 hours)<br>**Sun:** 4:00 to 22:00 (18 hours) | $r_1^{shoulder} = 0.107$ |
| Off-peak | Other hours | $r_1^{off-peak} = 0.060$ |

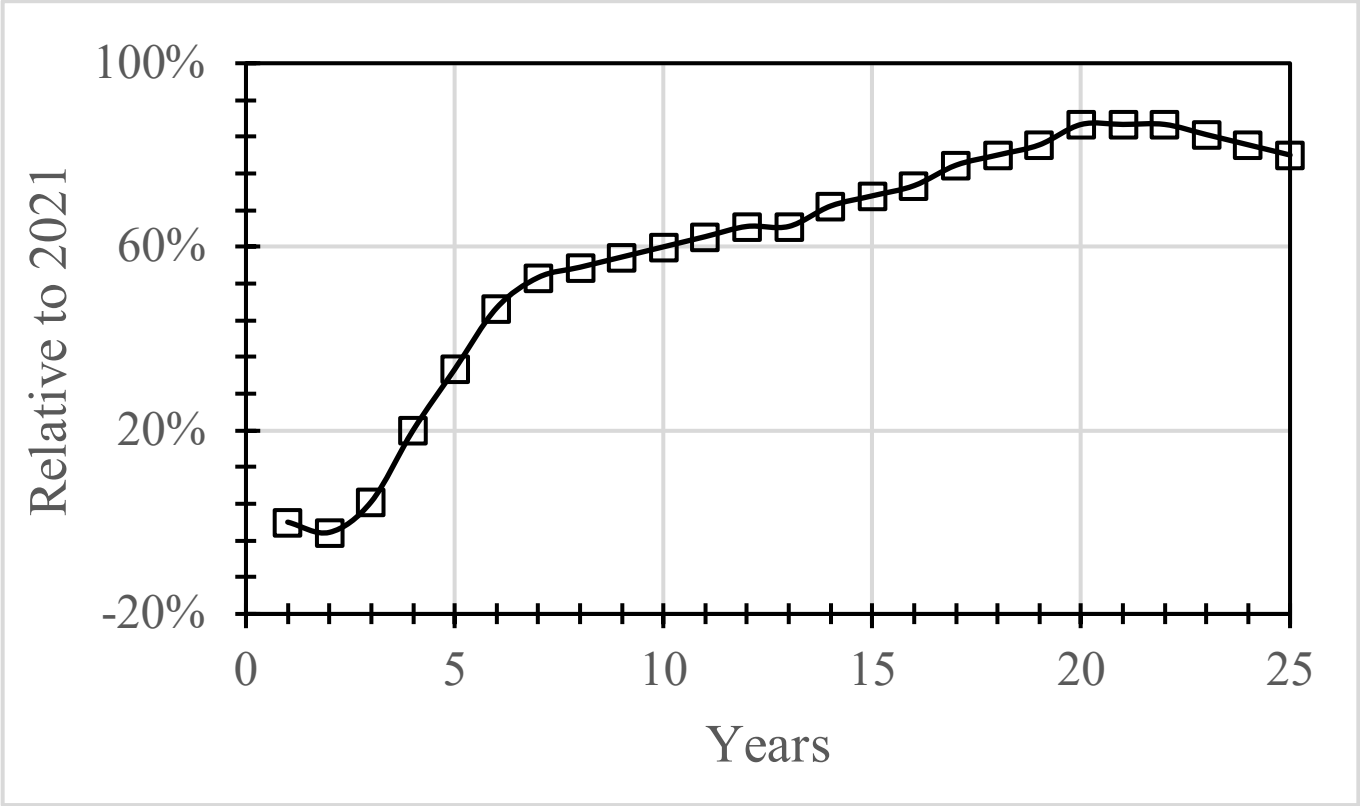

Figure 8 Electricity price projection over 25 years [65]

# 4 Result and discussion

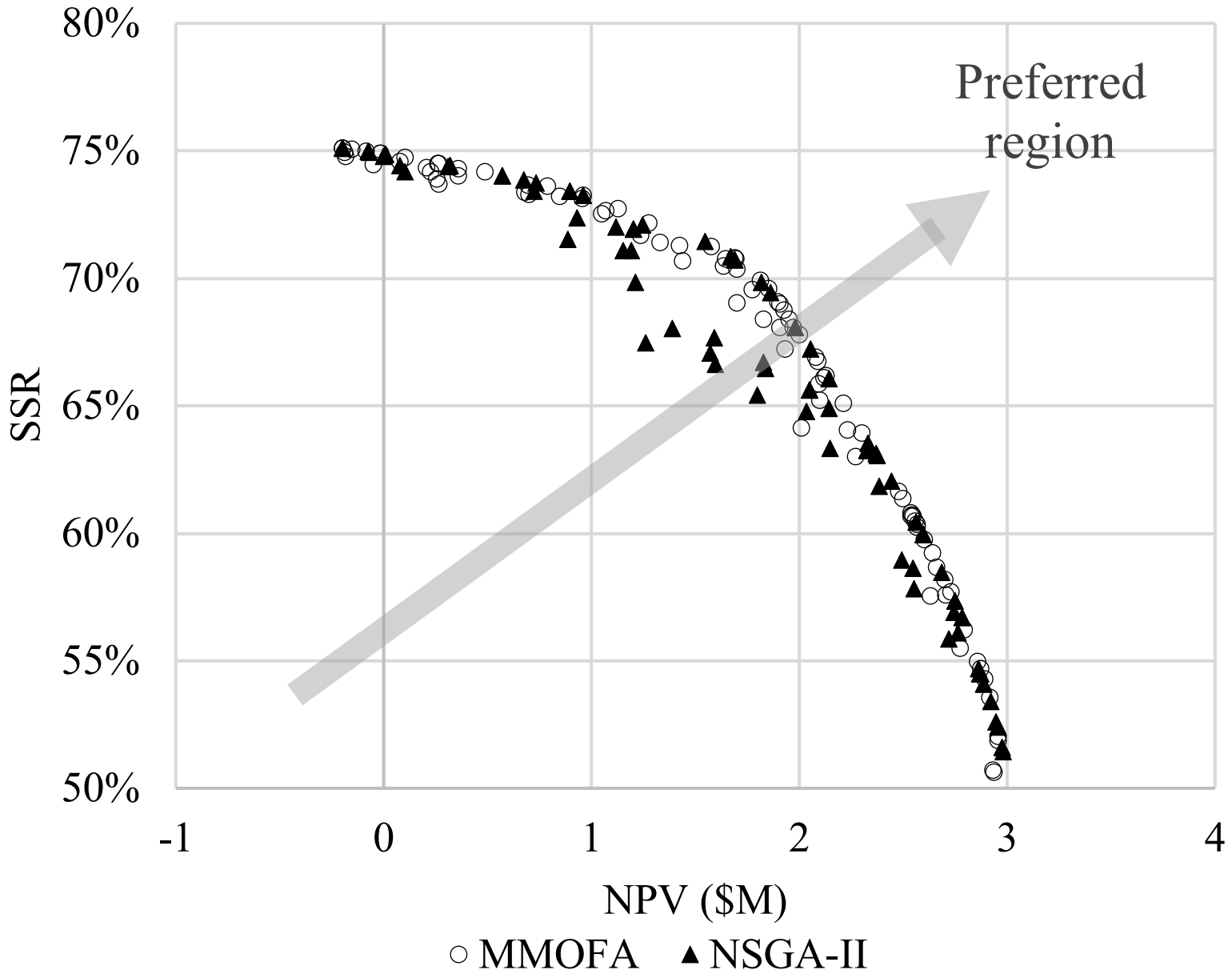

Figure 9 Superimposed results of five runs each of MOMFA and NSGA-II

In this section, the accuracy and performance of MMOFA proposed in this study are firstly compared with the NSGA-II algorithm. To accurately compare the performance of two optimisation algorithms, each algorithm was run five times to ensure that the results were not biased because of randomness in the simulation. Furthermore, the simulation was conducted on the same computer with the same computational time (2 hours for each run). Figure 9 shows the superimposed Pareto Front of five runs each of MOMFA and NSGA-II using Case 1C (see Figure 1). The performance of each algorithm is assessed primarily by the placement of its computed-optimal solutions closer to a preferred region with a higher NPV and a higher SSR - the closer the algorithm to the top right of the figure, the better performance. Figure 9 shows that MOMFA outperforms NSGA-II. While MOMFA has located and densely filled all parts of a single Pareto front, NSGA-II forms two separate fronts around the region between SSR =

65% and SSR = 72%. In order words, the optimisation results of NSGA-II are not consistent and generalised. This discrepancy is governed by the starting conditions, which include the initial random set of the algorithm. Therefore, there is a risk that NSGA-II has some runs that fall short of finding the local Pareto front, necessitating multiple runs of the algorithm to ensure it has identified the global Pareto front. Figure 8 shows that the MOMFA is a robust and accurate tool to optimise the ESS in this study. As a result, the MOMFA will be used to determine the optimal solutions for different ESS options (BESS/HESS with different operational strategies) in the following sections.

### 4.1 Energy storage option evaluation

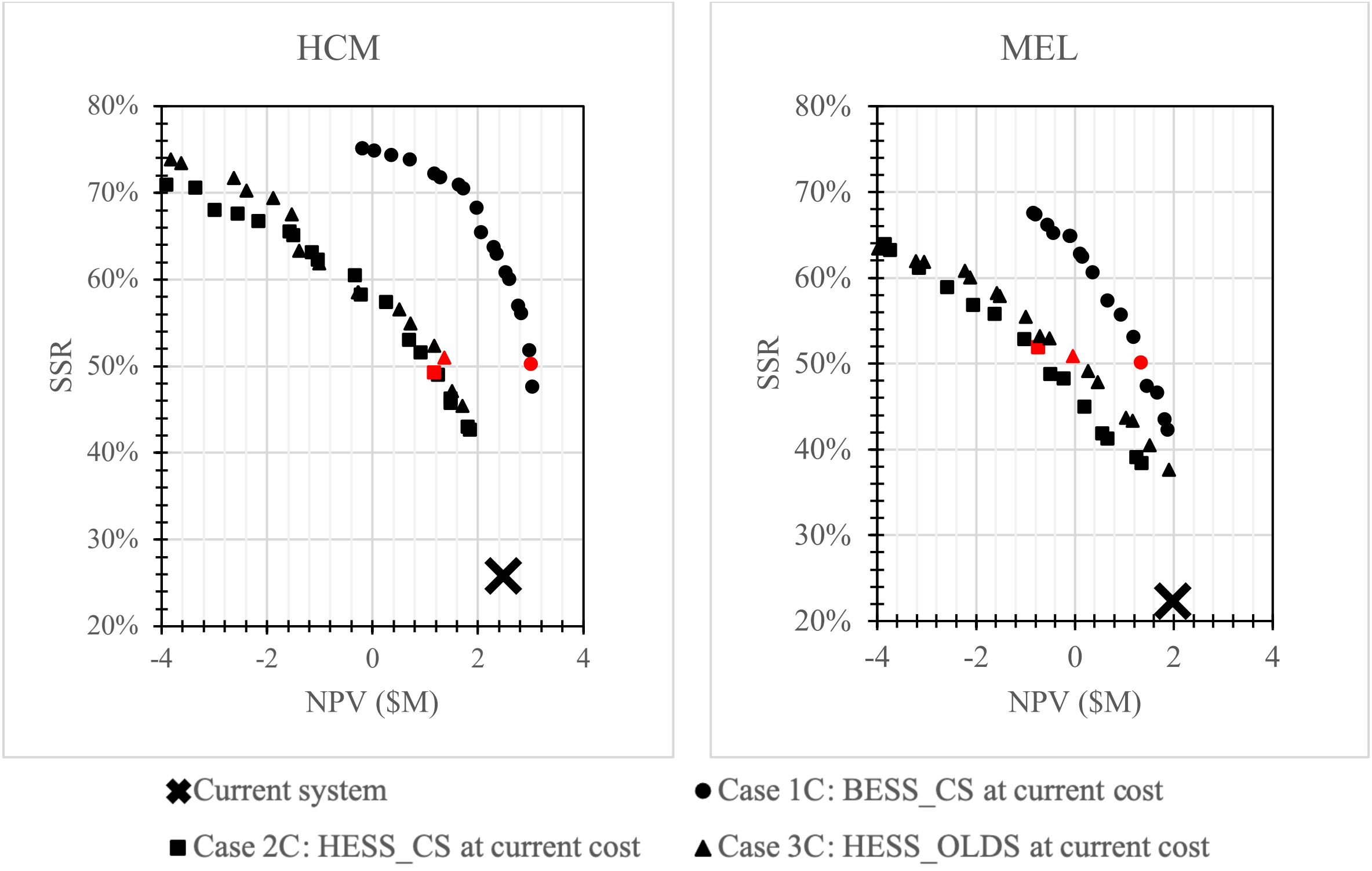

Figure 10 Pareto fronts of different storage options at current cost scenarios, red points indicate the systems are chosen in Table 5

This section presents the performance evaluation of the ESS options. The Pareto fronts of BESS at the current cost (Case 1C) are compared with those of HESS (Case 2C and 3C) in Figure 10. All options significantly increase SSRs in both locations. For the tropical location (i.e., HCM) the projected SSR, or the percentage of renewable energy supply over the total electricity demand, of the current system is 26% over 25 years. When ESSs are implemented, SSRs increase to between 43% and 75%. However, the results also show that only BESSs in HCM can bring economic benefits to the project. The projected NPV of the current system is $2.4M over 25 years. An extreme solution, a solution in one of two ends of the Pareto Front, of

the BESS system in HCM can increase the NPV up to $3.0M which is 22% higher than the current system. Meanwhile, extreme solutions of HESSs can only achieve an NPV around $1.7 to $2.0M. On the other hand, the projected NPV and SSR of the current systems are $2.0M and 22% over 25 years, respectively for the sub-tropical location (i.e., MEL). With the adoption of ESSs, SSRs span between 37% and 68%. The maximum NPV of extreme solutions in all three proposed options is only $1.9M. It is worth mentioning that all systems yield a lower NPV and SSR in sub-tropical weather compared to tropical weather. Figure 9 also shows that the BESS outperformance the HESS in the current cost scenario, indicating the cost challenges of hydrogen solution. However, the viability of HESS can be improved for locations with high seasonal variations. In fact, there is a significant separation between the optimal BESS and the two optimal HESS curves for HCM, which is closer for MEL.

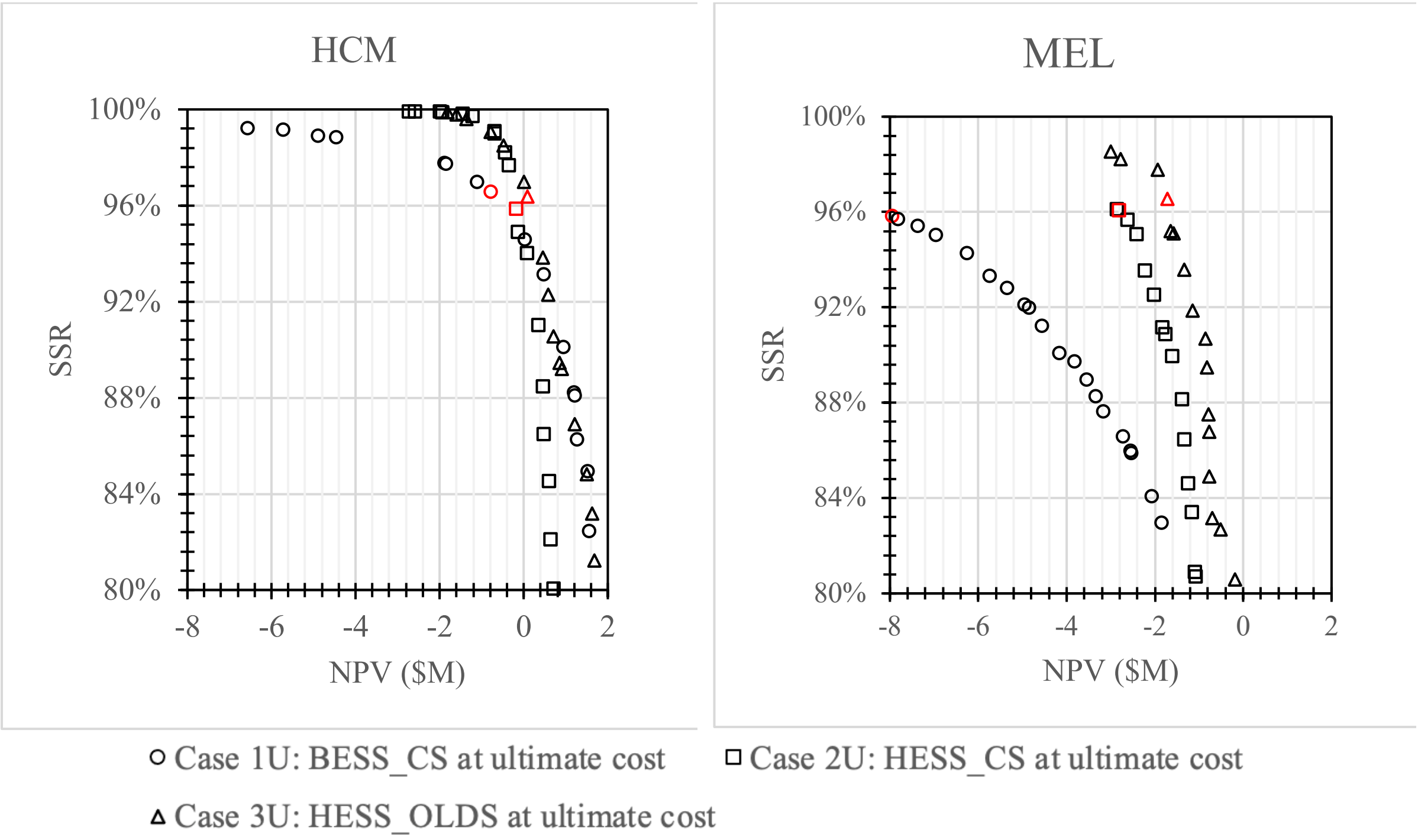

Figure 11 Pareto fronts of different storage options of high renewable energy penetration at ultimate cost scenarios, red systems are chosen in Table 5

The previous result (Figure 10) shows that the BESS, as a mature technology, gives a better performance than HESS with a higher NPV and SSR. However, this result is only applicable to energy storage that is needed in a short duration. With the high renewable energy penetration, this short time span is insufficient for reliable power system operation. Long-duration energy storage is required when renewable penetrations are above 80%, according to numerous studies focusing on high renewable power systems [66-69]. To investigate the role of long-duration energy storage in achieving the higher decarbonization of power systems in the future, a constraint is imposed to determine the optimal design with an SSR greater than 80%. The ultimate cost scenario is used in this analysis and the results are presented in Figure 11. It is interesting to observe that when the system approaches fully renewable (i.e., SSR > 95%), hydrogen becomes the preferable option in both locations. The results (Figure 11) indicate that there are no optimal solutions for the BESS that meets the full renewable energy requirement (i.e., SSR approaches 100%). In such requirements, the HESS would be a viable solution. Moreover, Figure 11 also confirms the advantages of HESS for highly seasonal locations. In Figure 11, in HCM, the BESS still performs well up to 95% while HESS dominates in MEL with an SSR of just 80%.

Table 5 Optimal configuration of different storage options with SRR around 50% and 95%

| | HCM | | | | | | MEL | | | | | |
|---|---|---|---|---|---|---|---|---|---|---|---|---|
| Case | 1C | 2C | 3C | 1U | 2U | 3U | 1C | 2C | 3C | 1U | 2U | 3U |
| PV (kWp) | 2,363 | 4,334 | 3,465 | 5,357 | 5,525 | 4,829 | 2,828 | 4,158 | 4,680 | 15,000 | 8,679 | 5,940 |
| Battery (kWh) | 1,890 | | | 11,072 | | | 3,480 | | | 9,015 | | |
| PEMEL (kW) | | 362 | 693 | | 4,320 | 2,581 | | 1,060 | 745 | | 2,755 | 3,270 |
| Storage Tank (kg) | | 90 | 374 | | 2,705 | 1,723 | | 254 | 411 | | 1,502 | 3,115 |
| PEMFC (kW) | | 67 | 88 | | 3,193 | 3,939 | | 206 | 105 | | 2,993 | 3,803 |
| $t_{start}$ (h d/m) | | | 3 AM<br>19/2 | | | 5 AM<br>4/8 | | | 3 AM<br>2/4 | | | 12 AM<br>31/3 |
| $t_{end}$ (h d/m) | | | 5 AM<br>31/3 | | | 10 AM<br>23/8 | | | 6 PM<br>15/9 | | | 4 AM<br>15/9 |
| LIMIT_SUNNY (% tank) | | | 100% | | | 100% | | | 26% | | | 16% |
| LIMIT_CLOUDY (% tank) | | | 24% | | | 30% | | | 100% | | | 100% |
| **NPC (M USD)** | 3.90 | 4.42 | 5.50 | 11.30 | 11.36 | 12.39 | 5.30 | 6.57 | 6.96 | 20.36 | 13.37 | 13.50 |
| **NPV (M USD)** | 3.00 | 1.16 | 1.36 | -0.78 | -0.18 | 0.07 | 1.33 | -0.75 | -0.05 | -7.96 | -2.85 | -1.72 |
| **SSR** | 50% | 49% | 51% | 96% | 96% | 96% | 50% | 52% | 51% | 96% | 95% | 95% |

To provide a detailed understanding of the optimal design of ESSs presented in Figure 10 and Figure 11, some optimal systems are selected to examine the sizing, overall cost and performance of each system. For the current cost scenario (Case 1C, 2C and 3C), each system is chosen with a similar SSR of around 50%. Similarly, for the ultimate cost scenario, the chosen systems have SSRs of around 95%. Table 5 depicts the optimised configuration of these systems. As seen in Table 5, the system sizes in the sub-tropical climate are always larger, and the NPCs are 10 to 80% greater than those systems in the tropical climate, under a same case and cost scenarios. Also, for both locations, the OLDS greatly increases the NPV of the HESS. As shown in Table 5, in HCM, the OLDS increases the NPV of the HESS by 17% from $1.16M to $1.36M. In MEL, the NPV of the system increases by 32% from $0.79M to $1.04M, highlighting the importance of optimising operational strategy. For the current cost scenario with an SSR of 50%, BESS always achieve a higher NPV than HESS due to a lower NPC. For a system with an SSR of 95%, it's not possible to achieve an economic solution because all storage technologies yield a negative NPV. Even though the cost of each component has been significantly decreased under the ultimate cost scenarios, high renewable energy penetration necessitates a considerably large component. Especially in the sub-tropical zone, when the solar radiation is low in cloudy months, the BESS needs a massive solar system of 15,000 kWp in order to meet the load demand. For HESS, the ability to store energy during sunny months and release it later during cloudy month alleviate the need for such a large solar PV system and only requires a system of 5,555 kWp. As a result, BESS requires a massive investment of $20.36M while HESS requires $13.50, leading to an NPV of -$1.72M for HESS and the NPV of -$7.96M for BESS. This highlights that HESS can be considered as a more viable solution in the future when the higher penetration of renewable energy is inevitable.

Figure 12 shows the breakdown NPC of selected systems in HCM in Table 5. For system with SSR of 50% (Case 1C, 2C, 3C), solar PV is the largest expenditure and usually accounts for more than half of NPC. Both PEMEL and PEMFC equally share a large portion while the storage tank accounts for only a minor fraction of the total HESS cost. For system with SSR of 95% (Case 1U, 2U, 3U), the share of the solar PV drops significantly as ESS becomes more prominent. Furthermore, the replacement cost of components accounts for a significant portion because replacement cost factor $\boldsymbol{r_{rep,i}}$ is not used in the ultimate cost scenario (due to the limits of raw material costs), indicating that future research should not only focus on lowering the cost, but also on increasing component longevity.

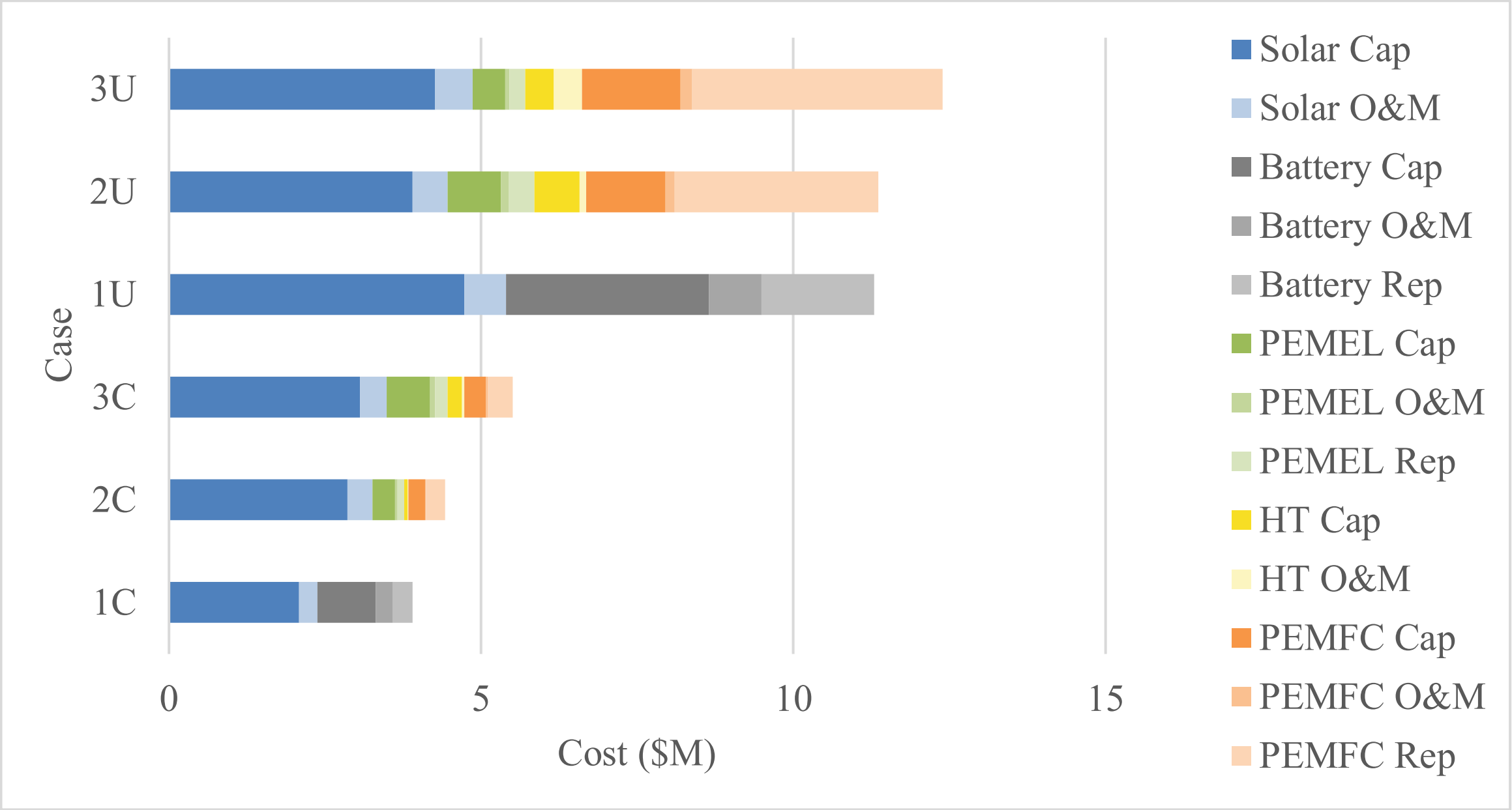

Figure 12 Break down NPC of selected systems in HCM. (For interpretation of the references to color in this figure legend, the reader is referred to the web version of this article.)

### 4.2 Operational strategy comparison

The results presented in the previous section demonstrate the superior performance of BESS for the current cost scenarios with low-to-medium SSR (i.e., <80%). However, the HESS with the long-term storage strategy has the potential to achieve higher penetration of renewable energy (i.e., SSR > 80%), which necessitates a deeper understanding of the impact of long-term storage of HESS. As a result, four systems from Case 2C and 3C in Table 5 are chosen to investigate the differences between the conventional and long-term operational strategies of HESS in HCM and MEL. Figure 13 presents the charge power of PEMEL, discharge power of PEMFC, and the storage level of the HESS in the first year for MEL location. MEL has a strong seasonal mismatch, which results in high solar output during the summer (Dec-Feb) and very low during winter (Jun-Aug) as depicted in Figure 13. This leads to a large variation in the charging power of PEMEL throughout the year as shown in Figure 13a. Similarly, the monthly discharge power of PEMFC is also lower in the OLDS, compared to the CS due to the controlled discharge operation of the OLDS. While the system with the CS will discharge whenever there is no excess solar, the OLDS stores hydrogen from sunny months and only uses it during peak and shoulder hours during cloudy months, indicating a cost-effective approach to dispatch the stored hydrogen. The level of storage in Figure 13c also shows that in CS, the storage tank level is low most of the time, especially between March and November, as this strategy just uses all available resources. The OLDS controls the release level and, therefore, maximizes the chance that there is hydrogen to discharge electricity during the peak and shoulder hours.

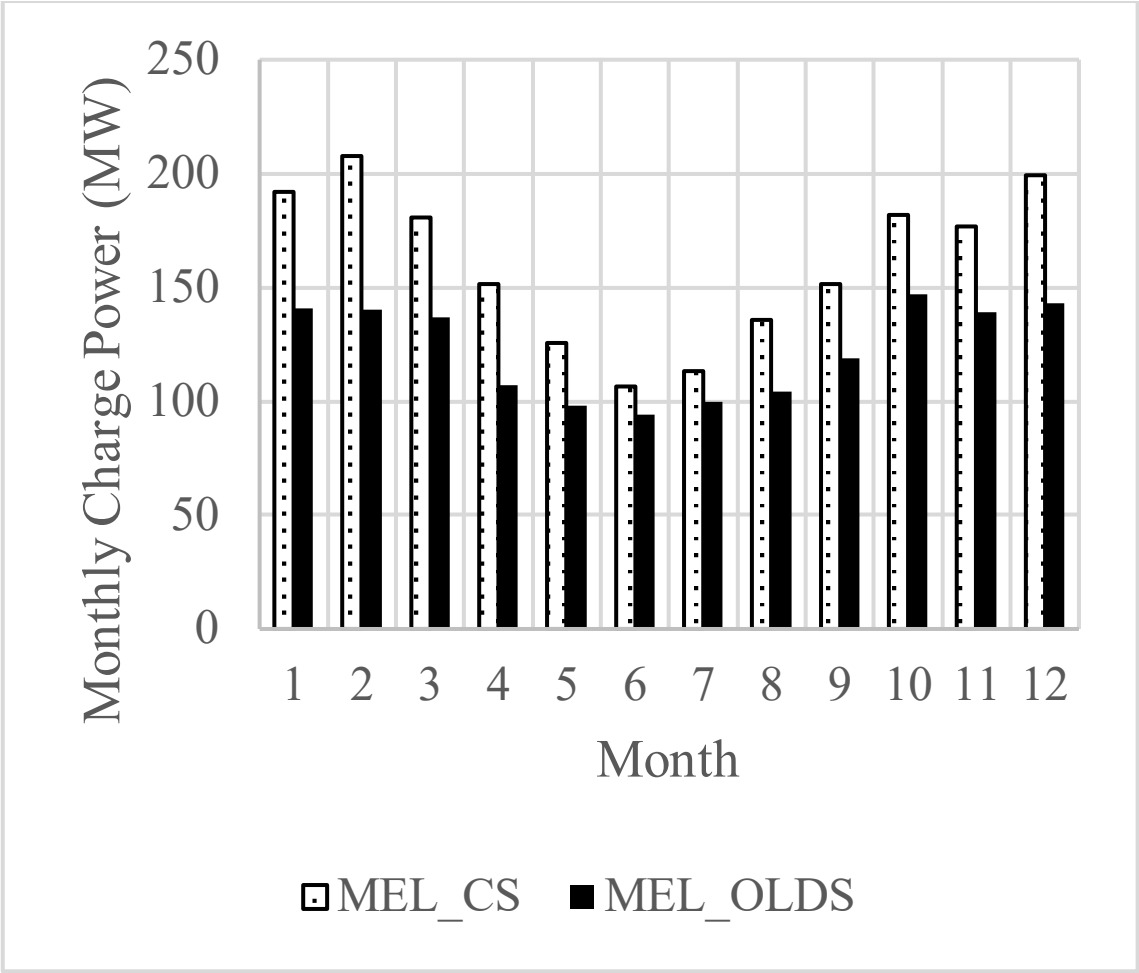

(a) Charge power of PEMEL

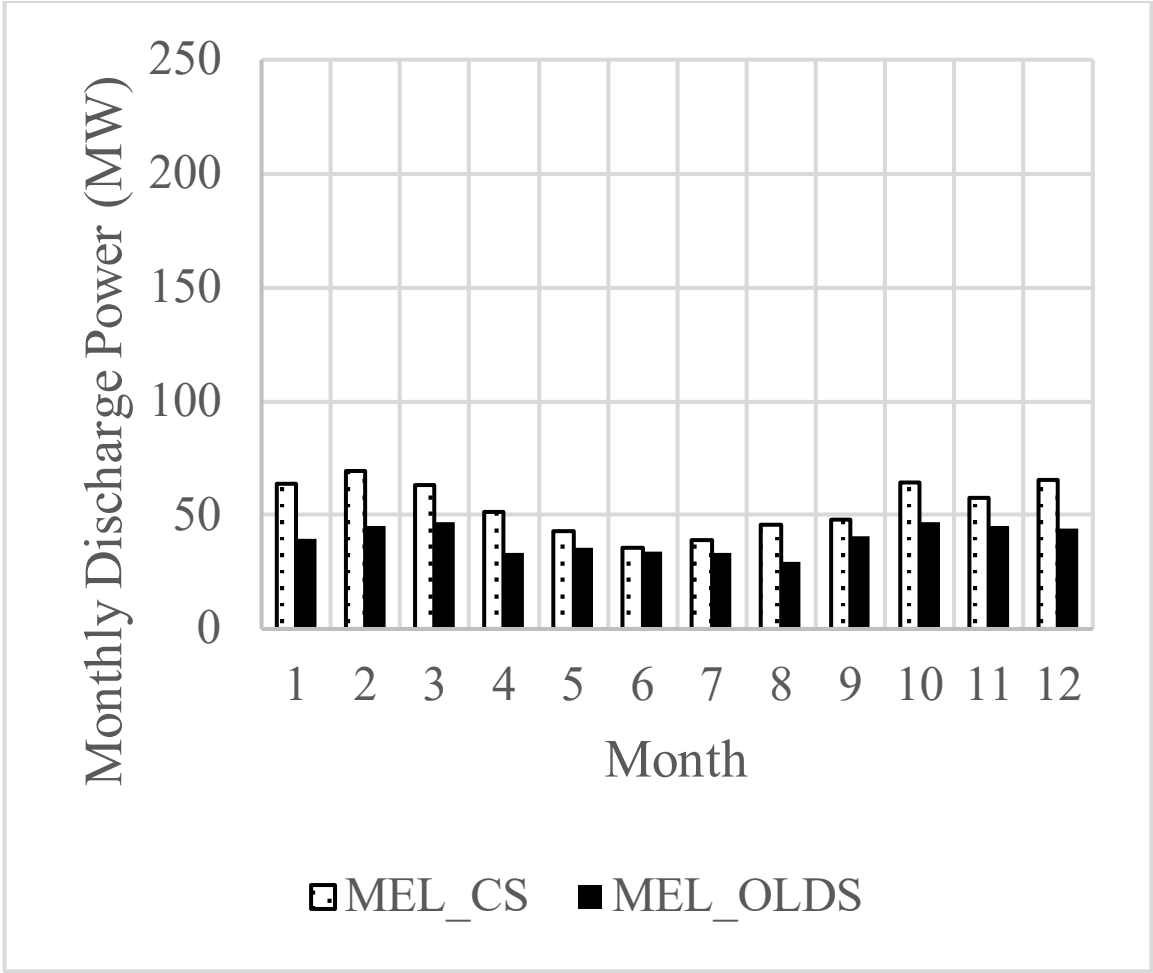

(b) Discharge power of PEMFC

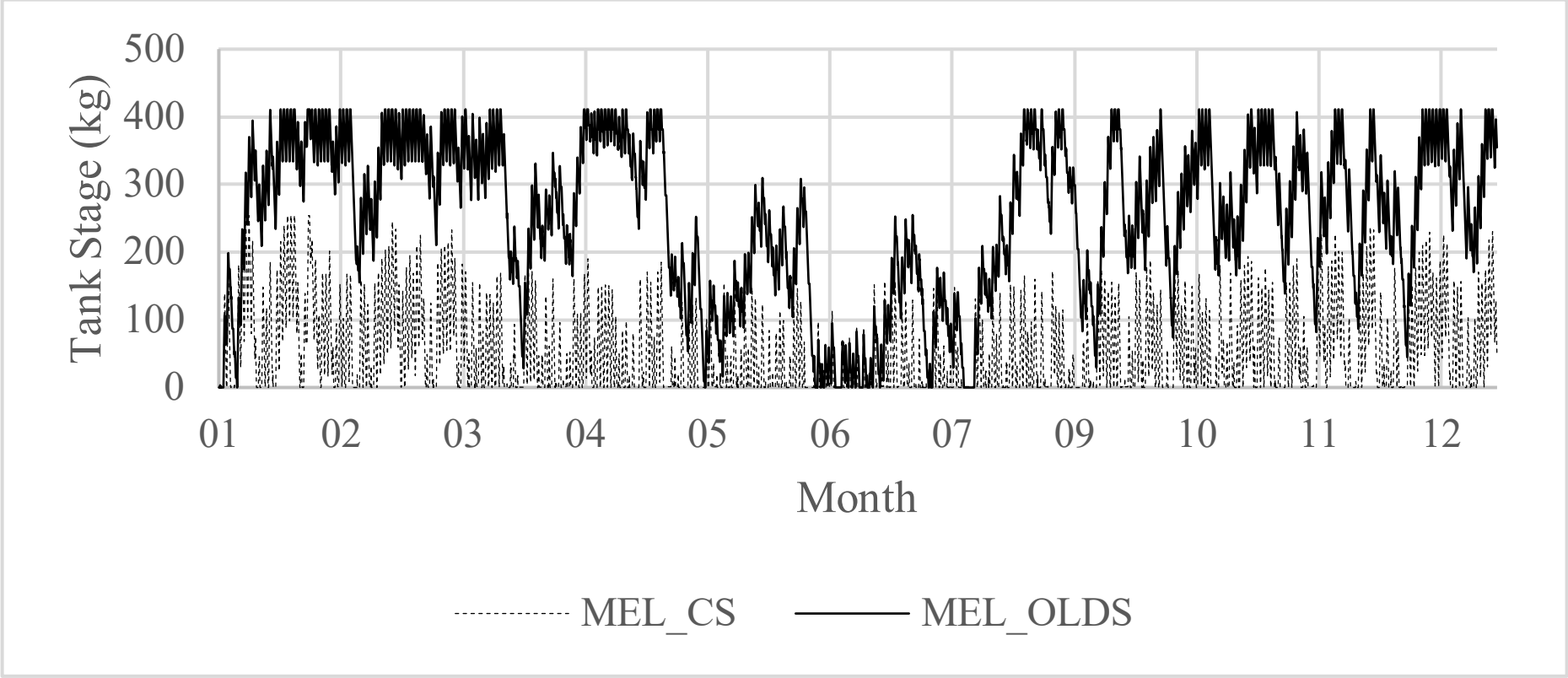

(c) Tank stage

Figure 13 Detailed operation of HESS for the MEL case in the first year

Technical performance from the first year of the simulation in HCM is presented in Figure 14. As HCM has year-round sunshine, the charging power from the PEMEL varies little throughout the year as shown in Figure 14a. In tropical climate zone, the two strategies address the optimisation problems differently as indicated by the size of the system shown in Table 5. In CS, the algorithm maximises the NPV by reducing NPC with a smaller system. The OLDS, on the other hand, aims to maximise the savings in electricity bill by ensuring that there is always enough hydrogen to generate electricity during peak and shoulder hours. As component sizes in OLDS are larger than those in CS, Figure 14a and b show that OLDS has a higher monthly charge and discharge power than CS. The amount of storage in Figure 14c further demonstrates that in CS, the storage tank level is frequently empty, whereas in OLDS, there is always ample storage.

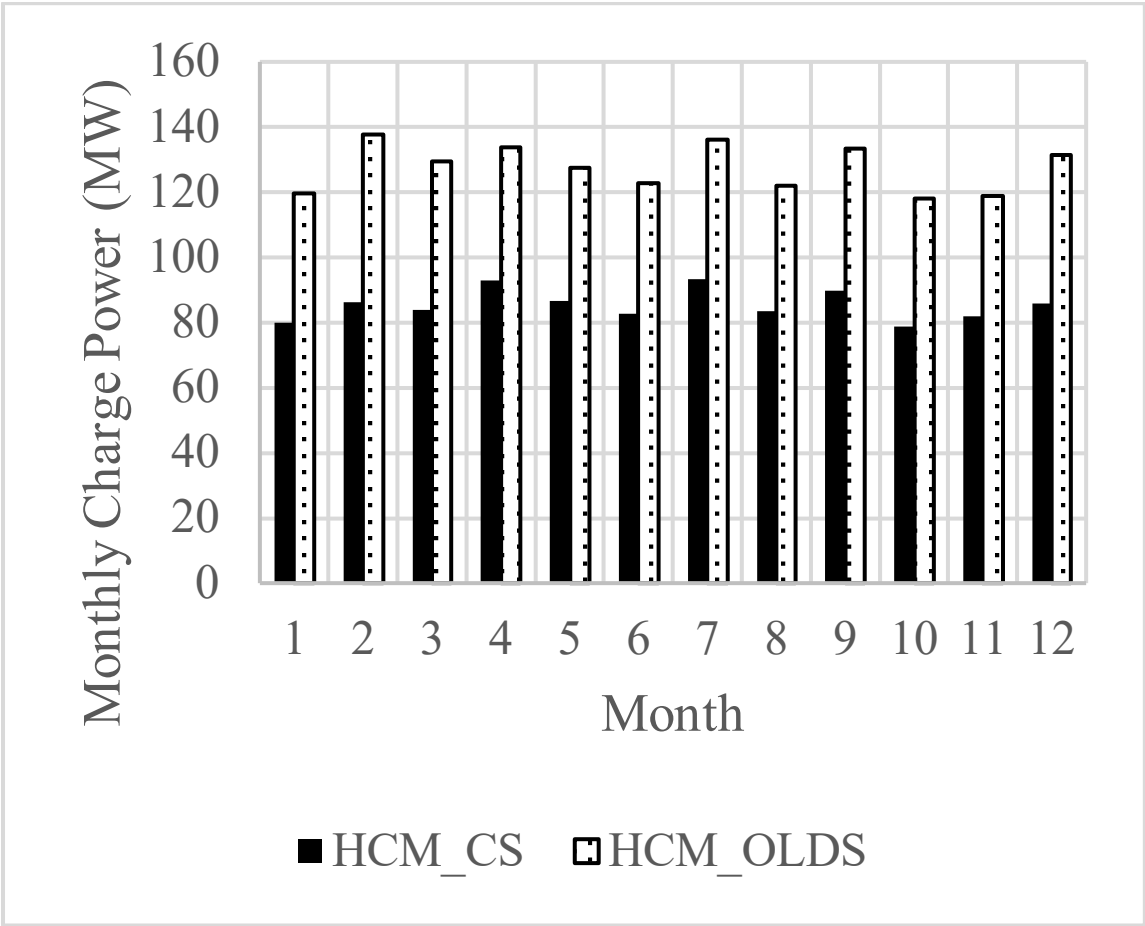

(a) Charge power of PEMEL

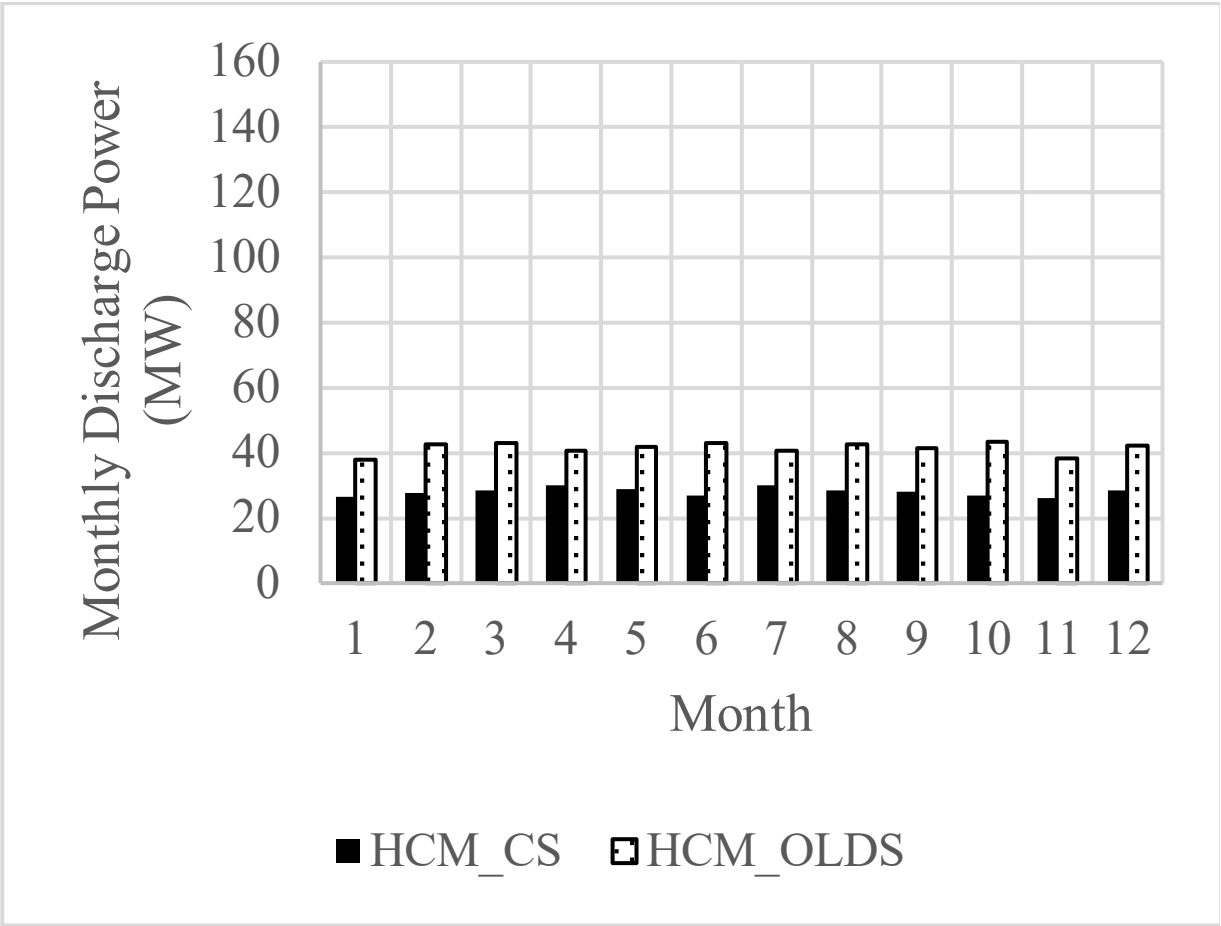

(b) Discharge power of PEMFC

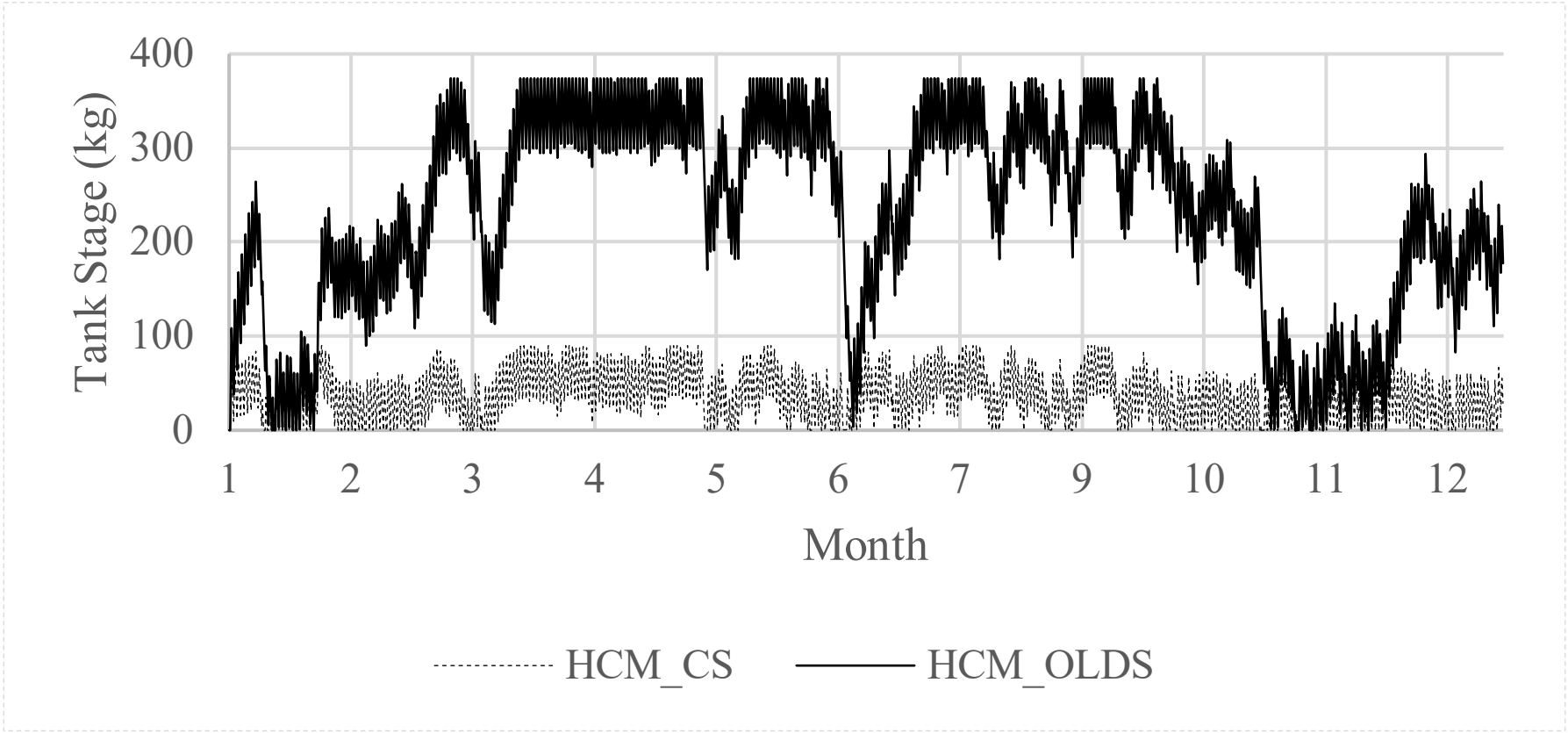

(c) Tank stage

Figure 14 Detailed operation of HESS for the HCM case in the first year

The degradation processes over 25 years of PEMEL and PEMFC are illustrated in Figure 16. Increases in SOH are due to replacements. At the replacement time, the SOH of both PEMFC and PEMEL is usually around 80 to 85%. There is no significant difference in PEMELs degradation between the two strategies in MEL. Due to the low solar radiation in MEL, the two PEMELs just operate at their potential. However, in HCM, the PEMEL with OLDS has a significant higher SOH than the system with CS. This happens as the storage in OLDS usually reached full capacity and cannot store more energy, resulting in less operation time than CS. In the case of PEMFC, the SOH of MEL_OLDS is higher than MEL_CS, indicating a better performance of OLDS. The SOH of HCM_OLDS is slightly lower than HCM_CS because the OLDS generates substantially more electricity than the CS, as shown by the discharge power in Figure 14c. In general, the OLDS effectively controls the operation of the HESS, reducing the degradation of system components and ultimately, increasing the total efficiency during the project lifetime.

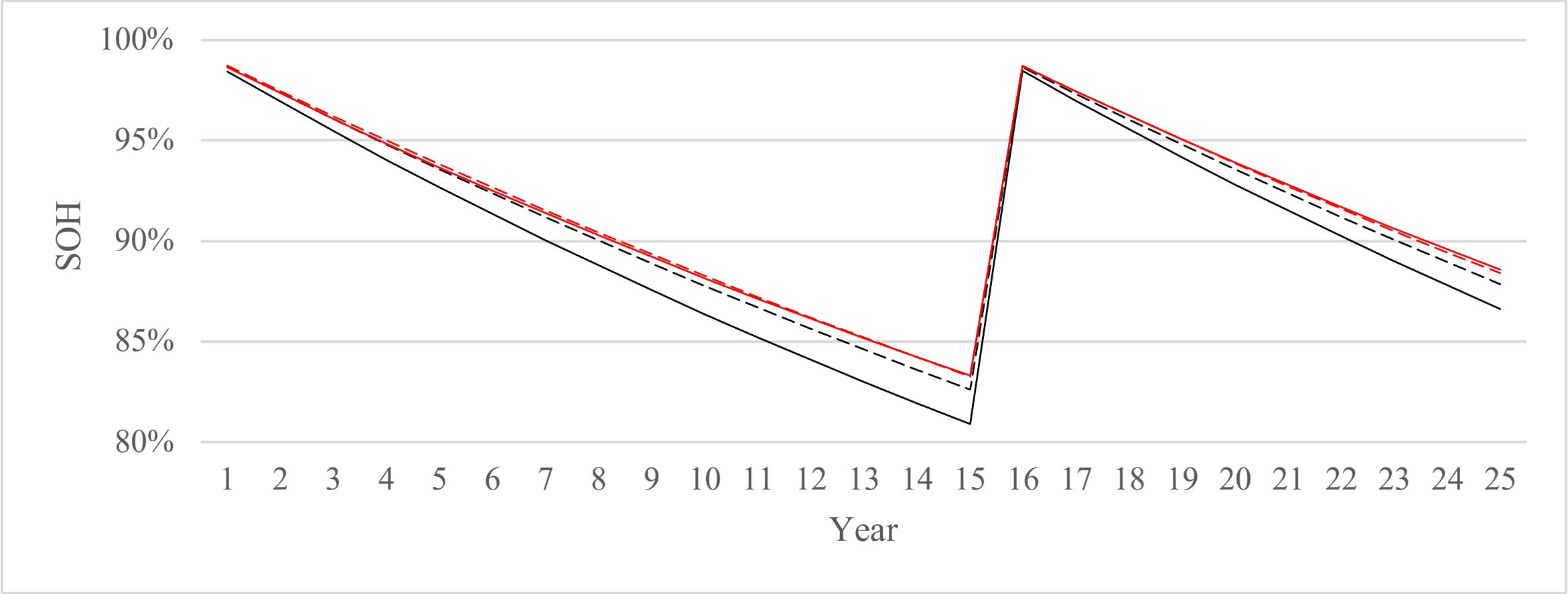

(a) Yearly degradation of PEMEL

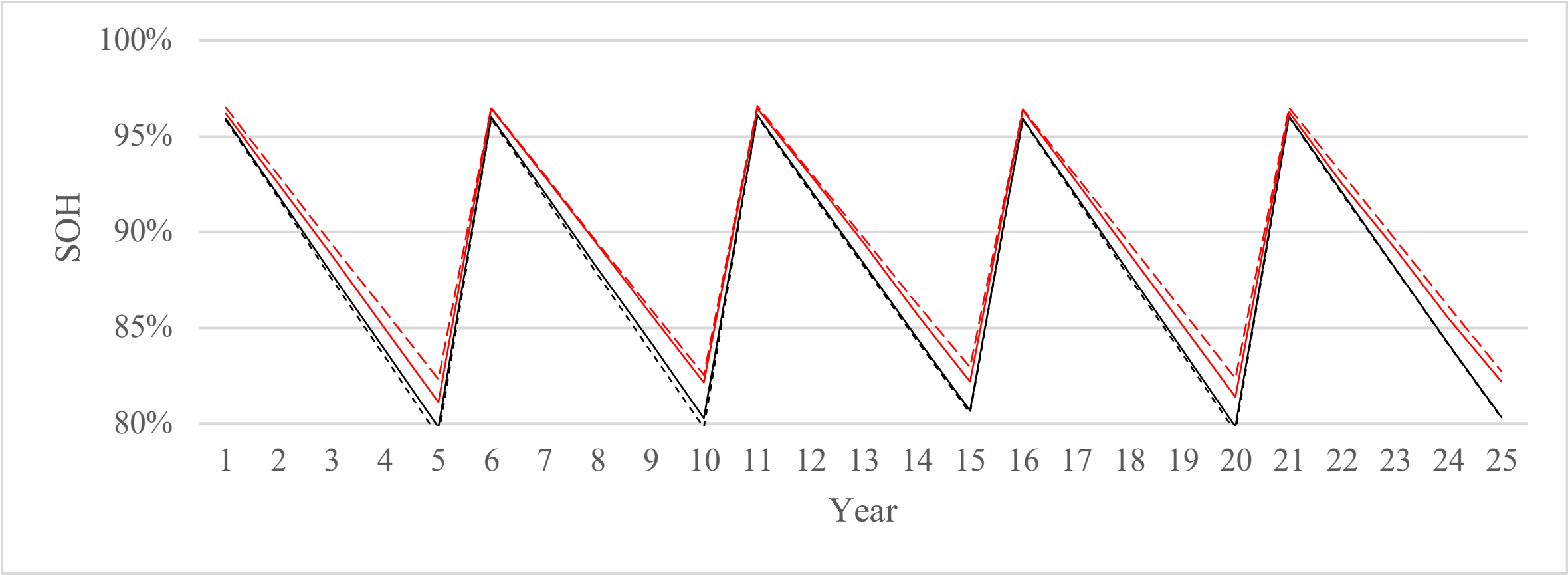

(b) Yearly degradation of PEMFC

HCM_CS HCM_OLDS MEL_CS MEL_OLDS

Figure 15 Evolution of the SOH of the PEMFC, PEMEL over 25 years.

## 5 Conclusion

This paper presents the optimisation study of sizing and operational strategy of a grid-connected PV-hydrogen/battery storage system using the Multi-Objective Modified Firefly Algorithm (MOMFA). The long-term storage potential of hydrogen energy storage system (HESS) is investigated by proposing the optimised long duration strategy (OLDS). A real-world case study in the tropical weather zone with actual, monitored data is used to examine the efficiency of the developed optimisation methodology with the consideration of electricity price forecasting, solar output depreciation, and degradation of storage components in the system. In addition, another location in the sub-tropical weather zone (with high seasonal variations) is hypothetically analysed to examine the performance of the system in different geographic locations.

The results show that the present MOMFA is more accurate and robust, compared with the popular NSGA-II algorithm. Comparing the two storage options: Battery (BESS) and Hydrogen energy storage system (HESS), it is found that the BESS outperforms the HESS with the current cost scenario and low-to-medium SSR. The solar system with the BESS yields a better Pareto front curve (higher SSR and NPV), compared with the system with the HESS. Nonetheless, the HESS is a viable option for the high requirement for using renewable energy (SSR>80) with the ultimate cost scenario, especially with using the OLDS for a geographic location with high seasonal variations. Further investigation on the operation of HESS with the proposed operational strategy reveals that the OLDS effectively controls the functioning of the HESS, lowering system component degradation and, as a result, enhancing the efficiency of the storage system during the project lifetime. The detailed analyses and results presented in this study can support researchers and engineers in developing storage systems for renewable energy.

## Acknowledgment

This research was supported by the Melbourne Postdoctoral Fellowship and the ARC Linkage Project (ARC LP210200774). The authors thank Copper Mountain Energy (CME), Viet Nam for the support to accomplish this research. The first author, Tay Son Le, would like to thank the University of Melbourne for offering the Melbourne Research Scholarship.